\documentclass[a4paper,11 pt]{amsart}
\usepackage{srollens-en}[2006/10/02]

\title[Nilmanifolds]{Geometry of nilmanifolds with left-invariant complex structure and deformations in the large}

\author{S\"onke Rollenske}

\address{Dr. S\"onke Rollenske\\Department of Mathematics\\
 Imperial College London \\
 SW7 2AZ London\\
 United Kingdom}
\email{s.rollenske@imperial.ac.uk}

\renewcommand{\lg}{\ensuremath{\gothg}}
\newcommand{\lh}{\gothh}
\newcommand{\einsnull}[1]{{{#1}^{1,0}}}
\newcommand{\nulleins}[1]{{{#1}^{0,1}}}

\newcommand{\defobj}[1]{\textbf{#1}}

\newcommand{\refb}[1]{{\upshape (\ref{#1})}}
\newcommand{\refenum}[1]{\textup{(\textit{#1})}}

\newcommand{\margincom}[1]{\marginpar{\quad}}
\begin{document}

\begin{abstract}
The relation between nilmanifolds with left-invariant complex structure and iterated principal holomorphic torus bundles is clarified and we give criteria under which deformations in the large are again of such type. As an application we obtain a fairly complete picture in complex dimension three.

\end{abstract}

\maketitle

\tableofcontents
\section*{Introduction}

A very general question in the theory of complex manifolds is the following: let $M$ be a compact, differentiable manifold and let
\[\kc:=\{ J\in \End (TM)\mid J^2=-id_{TM},J\text{ a complex structure}\}\]
be the space of all  complex structures on $M$.
What can we say about $\kc$ and its connected components?

If $\kc$ is non-empty consider a compact, complex manifold $X=(M,J)$.
The theory developed by Kodaira and Spencer in the 50's \cite{kod-sp58} and culminating in the theorem of Kuranishi \cite{kuranishi62} succeeds in giving a rather precise description of a slice of $\kc$ containing $J$ which is transversal to the orbit of  the natural action of $\mathrm{Diff}^+(M)$, called the Kuranishi slice.

While we have this powerful tool for the study of small deformations there is no general method available to study the connected components of $\kc$. 

From another point of view we say that two compact, complex manifolds $X$ and $X'$ are directly deformation equivalent $X{\sim}_\mathrm{def} X'$ if there exists an irreducible, flat family $\pi:\km\to \kb$ of compact, complex manifolds over a complex analytic space $\kb$ such that $X\isom\inverse\pi (b)$ and $X'\isom\inverse\pi (b')$ for some points $b, b'\in \kb$. The manifold $X$ is said to be a deformation in the large of $X'$ if both are in the same equivalence class with respect to the equivalence relation generated by ${\sim}_\mathrm{def}$, which is the case if and only if both are in the same connected component of $\kc$.

Even the seemingly natural fact that any deformation in the large of a complex torus is again a complex torus has been fully proved only in 2002 by Catanese  \cite{catanese02}. In \cite{catanese04} he studies more in general deformations in the large of principal holomorphic torus bundles, especially bundles of elliptic curves. This was the starting point for our research.

It turns out (see \cite{cat-fred06}) that the right context to generalise Catanese's results  is the theory of left-invariant complex structures on nilmanifolds, i.e., compact quotients of nilpotent real Lie-groups by discrete subgroups.

Nilmanifolds with left-invariant complex structure provide an important source for examples in complex
differential geometry.  Among these are the so-called Kodaira-Thurston manifolds, historically the first examples known to admit both a complex structure and a symplectic structure but no K\"ahler structure.
 In fact, a nilmanifold $M$ admits a K\"ahler structure if and only if 
it is a complex torus \cite{ben-gor88,hasegawa89} and the author showed \cite{rollenske07a}  that nilmanifolds can  be arbitrarily far from K\"ahler manifolds in the sense that the Fr\"olicher spectral sequence may be arbitrarily non-degenerate.


But unfortunately, even if every (iterated) principal holomorphic torus bundle can be regarded as a nilmanifold, the converse is far from true. Moreover it turns out that even a small deformation of a principal holomorphic torus bundle may not admit such a structure (see Example \ref{badex}).
Small deformations of nilmanifolds were already studied for so-called abelian complex structures in \cite{con-fin-poon06, mpps06} and for  general left-invariant complex structure on nilmanifolds in \cite{rollenske08c}. Thus we are concerned with these two problems:
\begin{itemize}
\item Give conditions under which a left-invariant complex structure on a nilmanifold gives rise to a structure of iterated (principal) holomorphic torus bundle.a
\item Study deformations in the large of iterated (principal) holomorphic torus bundles.
\end{itemize}

There is already a vast literature concerning nilpotent Lie-algebras and  left-invariant complex structures on nilmanifolds (see e.g. the articles of Console, Cordero, Fernandez, Fino, Grantcharov, Gray, McLaughlin, Pedersen,  Poon,  Salamon, Ugarte, et al. cited in the bibliography) and we recapitulate the basic results in Section \ref{LICS}, emphasising the complex geometric structure. In Section \ref{smalldefo} we will also recall the results on Dolbeault-cohomology and small deformations.

A nilmanifold $M$ can be characterised  by a  triple $(\lg, J, \Gamma)$ where $\lg$ is the nilpotent Lie-algebra associated to a simply connected nilpotent Lie-group $G$, $J$ is an integrable complex structure on $\lg$  and $\Gamma\subset G$ is a (cocompact) lattice. The datum of either $\lg$ or $\Gamma$ (considered as an abstract group) determines $G$ up to unique isomorphism. 
The general philosophy is that the geometry of the compact, complex manifold $M_J=(\Gamma\backslash G,J)$ should be completely determined by  the linear algebra of $\lg$, $J$ and the $\IQ$-subalgebra generated by $\log \Gamma\subset \lg$. By abuse of notation we will sometimes write $M_J=(\lg, J, \Gamma)$.

The real geometry of nilmanifolds is well understood but only the existence of  what we call a stable  torus bundle series in the Lie-algebra $\lg$ (Definition \ref{stableseries}) gives us sufficient control over the complex geometry of $(M,J)$. Geometrically this notion means the following:

On any real nilpotent Lie-group $G$ there is a filtration by normal subgroups, e.g., the ascending central series, 
\[G  \supset H_s\supset\dots\supset  H_1\supset H_0=\{1\}\]
such that $H_k/H_{k-1}$ is abelian and $\Gamma\cap H_k$ is a lattice in $H_k$ for any lattice $\Gamma\subset G$.

In other words, any compact nilmanifold $M=\Gamma\backslash G$ can be represented as a tower of differentiable  torus bundles $\pi_k: M_k\to M_{k+1}$
\[\xymatrix{ T_1 \ar@{^(->}[r] & M_1=M \ar[d]^{\pi_1}\\
T_2 \ar@{^(->}[r] & M_2\ar[d]^{\pi_2}\\
&\vdots\ar[d]\\
T_{s-1}\ar@{^(->}[r] & M_{s-1}\ar[d]^{\pi_{s-1}}\\
& M_s=T_s}\]
where the compact torus $T_k$ is the quotient of $H_k/H_{k-1}$ by the lattice and $M_k$ is the compact nilmanifold obtained from the Lie-group $G/H_{k-1}$ by taking the quotient with respect to the image of $\Gamma$.

Such a decomposition into an iterated principal bundle is far from unique. If the subgroups can be chosen in such a way that for every left-invariant complex structure $J$ on $M$ we have
\begin{enumerate}
\item The complex structure $J$ induces a left-invariant complex structure on $M_k$ for all $k$.
\item All the maps $\pi_k:M_k\to M_{k+1}$ are holomorphic with respect to these complex structures.
\end{enumerate}
 then we say that $\lg$ admits a stable torus bundle series. 


This phenomenon is not as uncommon as it may seem at first sight; for example every principal holomorphic torus bundle over an elliptic curve which is not a product has this property (Example \ref{torusovercurvedescription}). Section \ref{structureclass} will be devoted to  giving a fairly complete picture of the occurring cases if $\dim [\lg,\lg]\leq 3$.
A large part of the classification of complex structures on real 6-dimension nilpotent Lie-algebras  \cite{salamon01, math.DG/0411254} can be recovered from our more general results.

It turns out that finding a stable torus bundle series (STBS) for some nilmanifold $M_J$ is a good step on the way to prove that every deformation in the large of $M_J$ is again such a nilmanifold. Indeed, in this case the holomorphic fibration over a torus $\pi: M\to T_s$ can be realised as a (topologically) fixed quotient of the Albanese variety and this will enable us in Section \ref{albanesequotients} to determine all deformations in the large if the fibres of $\pi$ have sufficiently nice properties:

\begin{custom}[Theorem \ref{nillargedeform}] 
Let $G$ be a simply connected nilpotent Lie-group with Lie-algebra $\lg$ and let $\Gamma\subset G$ be a lattice such that the following holds:
\begin{enumerate}
\item $\lg$ admits a stable torus bundle series  $(\ks^i\lg)_{i=0,\dots, t}$ (cf. Definition \ref{stableseries}).
\item Nilmanifolds with Lie algebra $\ks^{t-1}\lg$ and fundamental group $\Gamma\cap \exp(\ks^{t-1}\lg)$ constitute a good fibre class (cf. Definition \ref{goodfibreclass}).
\end{enumerate}
Then any deformation in the large of a nilmanifold with left-invariant complex structure $M=(\Gamma\backslash G, J)$ is again of the same type. 
\end{custom}
Generalising the methods used in \cite{catanese04} we will have to deal with the fact  that the dimension of the Albanese variety may vary in a family of nilmanifolds. 

In Section \ref{applications} we apply our results on the deformation of complex structures to the classification  obtained in Section \ref{structureclass} and thus give several series of examples (in arbitrary dimension) for which left-invariant complex structures yield a union of connected components of the space of complex structures $\kc$. This holds for example if  $\dim [\lg,\lg] \leq 1$ or if  $\lg$ is 3-step nilpotent and $\dim [\lg,\lg]=2$. We refer to Theorem \ref{apptheosmall} and Theorem \ref{apptheo} for a complete list of cases.

There are only finitely many  real 6-dimensional, nilpotent Lie-algebras that admit a complex structure \cite{magnin86, salamon01}  and we obtain a quite complete picture of their geometry and deformations in Section \ref{dim3}:
\begin{custom}[Theorem \ref{3dim-class}]
Let $M=(\lg , J, \Gamma)$ be a complex 3-dimensional nilmanifold with left-invariant complex structure. 
\begin{enumerate}
\item If $\lg$ is not in $\{\lh_7, \lh_{19}^-, \lh_{26}^+\}$, then $\lg$ admits a SPTBS and hence $M_J$ has the structure of an iterated principal holomorphic torus bundle. We list the possibilities in the following table:
\begin{center}
\begin{tabular}{c|c|c}
 base& fibre & corresponding Lie-algebras\\
\hline
3-torus & - & $\lh_1$\\
2-torus & elliptic curve & $\lh_2, \lh_3  , \lh_4, \lh_5, \lh_6$\\
elliptic curve & 2-torus & $\lh_8$\\
Kodaira surface & elliptic curve & $\lh_9, \lh_{10}, \lh_{11}, \lh_{12}, \lh_{13}, \lh_{14}, \lh_{15}, \lh_{16}$\\
\end{tabular}
\end{center}
Every deformation in the large is of the same type.
\item 
If   $\lg=\lh_{19}^-$ or $\lg=\lh_{26}^+$  then $\lg$ admits a STBS and $M_J$ can be described as 2-torus bundle over an elliptic curve but there is no principal torus bundle structure. Every deformation in the large is of the same type.
\item If $\lg=\lh_7$ then there is a dense subset of the space of all left-invariant complex structures for which $M$ admits the structure of principal holomorphic bundle of elliptic curves over a Kodaira surface but this is not true for all complex structures.
\end{enumerate}
\end{custom}

It would also be interesting to study the space of left-invariant complex structures  more in detail, for example determine when it is smooth or universal,  and perhaps look for some kind of moduli space or a description of a connected component of  the Teichm\"uller space. These questions have already been addressed in several cases by various authors \cite{catanese04, cat-fred06, GMPP04, ket-sal04}.

Nevertheless the conditions given in  Theorem \ref{largedeform} are rather strong and we show in Example \ref{badex} that they do not need to be satisfied. This leads to the question:
Which is the simplest example of a nilmanifold $M_J$ with left-invariant complex structure such that not every deformation in the large carries a left-invariant complex structure?

Our theory of Albanese-Quotients in Section \ref{albanesequotients} is not restricted to nilmanifolds and it would be nice to find other applications.

\subsubsection*{Acknowledgements.} This work is part of my PhD-thesis \cite{rollenske07}. I would like to express my gratitude to my adviser Fabrizio Catanese for suggesting this research, constant encouragement and several useful and critical remarks. Simon Salamon gave  valuable references to the literature and helped to improve the presentation at several points. J\"org Winklemann remarked that every good Albanese variety is in fact very good. An invaluable bibliographic hint due to Oliver Goertsches opened a new perspective on the problem. The referee and the editor made several suggestion which helped to improve the presentation of the results.

I am grateful for the hospitality of the Max Planck Institut in Bonn where Section \ref{albanesequotients} got the final polishing; during the revision of the article I was visiting Imperial College London supported by a DFG Forschungsstipendium.

\section{Nilpotent Lie-algebras and nilmanifolds with left-invariant complex structure}\label{LICS}

In this section we will introduce the objects of our study and describe their basic properties. We will emphasise the complex-geometric structure of the nilmanifolds but the expert will find nothing new.

The geometrically important notion of stable (principal) holomorphic torus bundle series (SPTBS) will be given in Definition \ref{stableseries}.

\subsection{Lie-algebras with a complex structure}\label{basicdefin}
We will throughout need the yoga of almost complex structures and will now recall some basic definitions and notations.

Let $\lg$ be a finite dimensional real Lie-algebra and $J$ an almost complex structure on the underlying real vector space, i.e., $J$ is an endomorphism of $\lg$ such that $J^2=-Id_\lg$; setting $ix=Jx$ for $x\in \lg$ this makes $\lg$ into a complex vector space.
We can decompose the complexified Lie-algebra $\lg_\IC={\gothg}^{1,0}\oplus \gothg^{0,1}$ into the $\pm  i$-eigenspaces of the $\IC$-linear extension of $J$ and every decomposition $\lg_\IC=U\oplus \bar U$ gives rise to a unique almost complex structure $J$ such that  $\einsnull{\lg}=U$.

%
%
 
Usually we will use small letters $x,y,\dots$ for elements of $\lg$ and capital letters $X,Y,\dots$ for elements in $\einsnull\lg$; elements in $\nulleins\lg$ will be denoted by $\bar X,\bar Y,\dots$ using complex conjugation. 

The exterior algebra of the dual vector space $\lg^*$ decomposes as
\[\Lambda^k\lg^*_\IC=\bigoplus_{p+q=k}\Lambda^p{\gothg^*}^{1,0}\tensor \Lambda^q{\gothg^*}^{0,1}=\bigoplus_{p+q=k}\Lambda^{p,q}{\gothg^*}\]
and we have $\overline{\Lambda^{p,q}{\gothg^*}}=\Lambda^{q,p}{\gothg^*}$.
A general reference for the linear algebra coming with a complex structure is  \cite{Huybrechts} (Section 1.2).

\begin{defin}\label{defintegrable}
An almost complex structure $J$ on a real Lie-algebra $\lg$ is said to be \defobj{integrable} if the Nijenhuis condition
\begin{equation}\label{nijenhuis}
 [x,y]-[Jx,Jy]+J[Jx,y]+J[x,Jy]=0
\end{equation}
holds for all $x,y\in \lg$ and in this case we call the pair $(\lg , J) $  a Lie-algebra with complex structure. 
\end{defin}

Hence by a complex structure on a Lie-algebra we will always mean an integrable one; otherwise we will speak of almost complex structures. 
It is easy to show that $J$ is integrable if and only if $\lg^{1,0}$ is a subalgebra of $\lg_\IC$ with the induced bracket.

If $G$ is a real Lie-group with Lie-algebra $\lg$ then giving a left-invariant almost complex structure on $G$ is equivalent to giving an almost complex structure $J$ on $\lg$ and $J$ is integrable if and only if it is integrable as an almost complex structure on $G$. It then induces a complex structure on $G$ by the Newlander-Nirenberg theorem (\cite{Kob-NumII}, p.145) and $G$ becomes a complex manifold. The elements of $G$ act holomorphically by multiplication on the left but $G$ is not a complex Lie-group in general.

\subsection{Nilmanifolds with left-invariant complex structure }\label{set-up}

In this section we will collect a bunch of results on nilmanifolds, most of them well know and for which we claim no originality; but our presentation will emphasise the geometric structure of compact nilmanifolds.

If not otherwise stated $(\lg, J)$ will always be a Lie-algebra with (integrable) complex structure and $G$ will be a associated simply connected Lie-group. By a torus we will always mean a compact torus.

\begin{defin}
A compact complex manifold $M_J:=(M,J)$ is called \defobj{nilmanifold with left-invariant complex structure} if there is a nilpotent Lie algebra with complex structure $(\lg,J)$ and a lattice $\Gamma$ in an associated simply-connected Lie-group such that $M_J\isom(\Gamma\backslash G, J)$.
\end{defin}
Since $\Gamma=\pi_1(M_J)$ determines $G$ and hence $\lg$ up to isomorphism (\cite{VinGorbShvart}, p.45, Corollary 2.8\margincom{check the number}),  we will always identify  $M_J$ with $(\Gamma\backslash G, J)$ and call it a nilmanifolds with left-invariant complex structure of type $(\lg, \Gamma)$.
By abuse of notation we will sometimes write $M_J=(\lg, J, \Gamma)$.


For nilpotent Lie-groups the exponential map $\exp: \lg \to G$ is a diffeomorphism  and all analytic subgroups are closed and simply connected as well (\cite{Varadarajan}, Theorem. 3.6.2, p. 196).
The following  often gives the possibility to use inductive arguments.
\begin{lem}\label{fibration}
Let $(\lg, J)$ be a nilpotent Lie-algebra with complex structure. Let $\lh$ be an ideal in $\lg$ such that $J\lh=\lh$. Let $G$ and $H$ be the associated simply connected Lie-groups endowed with the left-invariant complex structures induced by $J$. Then there is a holomorphic fibration $\pi:G\to G/H$ with fibre $H$.
\end{lem}
\pf The map $\pi:G\to G/H$ is a real analytic fibration by the theory of Lie-groups and Lie-algebras since $H$ is closed. Hence it remains to show that the differential of $\pi$ is $\IC$-linear and since the complex structure is left-invariant it suffices to do so at the identity. But in this point the differential is given by the quotient map $\lg \to \lg/\lh$ which is $\IC$-linear by assumption.\qed

Remark that we used the nilpotency of $\lg$ only to ensure that $H$ is a closed subgroup.

\subsubsection{The real structure of $\Gamma\backslash G$}
We leave aside the complex structure for a moment and describe the geometry of the underlying real manifold.
In the Lie-algebra $\lg$  we have the following filtrations:
\begin{itemize}
\item The descending central series (nilpotent series) is given by 
\[ \kc^0\lg :=\lg, \qquad \kc^{i+1}\lg := [\kc^{i}\lg, \lg]\]
\item The ascending central series is given by 
\[\kz^0\lg:= 0, \qquad \kz^{i+1} \lg := \{ x\in \lg \mid [x,\lg ]\subset  \kz^{i}\lg\}.\]
In particular $\kz^1\lg=\kz\lg$ is the centre of $\lg$.
\end{itemize}
The Lie-algebra $\lg$  is called \emph{$s$-step nilpotent} if $\kc^s\lg=0$ and  $\kc^{s-1}\lg\neq 0$ or equivalently $\kz^s\lg= \lg$ and $\kz^{s-1}\lg \subsetneq \lg$.

Let $G_i$ be the simply connected nilpotent Lie-group associated to the Lie-algebra $\lg/\kz^{i-1}\lg$ and let $F_i\isom \IR^{n_i}$ be the abelian Lie-group corresponding to $\kz^{i}\lg/\kz^{i-1}\lg$. Then the map $G_i\to G_{i+1} $ is a principal bundle with fibre $F_i$ and, if $\lg$ is $s$-step nilpotent, we get a tower of such bundles:
\[\xymatrix{ F_1 \ar@{^(->}[r] & G_1=G \ar[d]\\
F_2 \ar@{^(->}[r] & G_2\ar[d]\\
&\vdots\ar[d]\\
F_{s-1}\ar@{^(->}[r] & G_{s-1}\ar[d]\\
& G_s=F_s.}\]
A manifold which admits such a tower of (principal) bundles is called a \defobj{$s$-step iterated (principal) bundle}.
We say that the structure of iterated (principal) bundle is induced by a filtration $\ks^k\lg$ in $\lg$ if we have a tower of (principal) bundles as above such that $F_i$ is (the quotient by the lattice of) the Lie-group corresponding to $\ks^{i}\lg/\ks^{i-1}\lg$ and 
$G_i$ is (the quotient by the lattice of the)  Lie-group associated to the Lie-algebra $\lg/\ks^{i-1}\lg$.

\subsubsection{The complex geometry of the universal covering G} 
To study the complex geometry of $(\lg, J)$ we need to study the interplay of $\lg$ and $J$. The descending central series associated to $J$  is defined by
\[\kc_J^i\lg:= \kc^{i}\lg+J\kc^{i}\lg \]
and it is known to have the following properties:
\begin{enumerate}
\item $\kc_J^i\lg$ is a $J$-invariant subalgebra of $\lg$ and an ideal in $\kc_J^{i-1}\lg$ \cite{con-fin01}.
\item We have always $\kc^i\lg \subset\kc_J^i\lg$ but in general inclusion can be strict. Nevertheless we have always $\kc_J^{1}\lg\neq \lg$. (\cite{salamon01}, Theorem. 1.3).
\item There is always a holomorphic fibration of $G$ over the abelian Lie-group (vector space) $\lg/\kc_J^1\lg$ whose typical fibre is the simply connected Lie-group $H$ associated to $\kc_J^1\lg$ with the left-invariant complex structure induced by the restriction of $J$. This is a real principal $H$-bundle but in general $H$ will not be a complex Lie-group so there is no way to speak of a holomorphic principal bundle.
One can iterate this process to find an iterated holomorphic bundle structure on $(G,J)$.
\end{enumerate}

The ascending central series associated to $J$ (minimal torus bundle series) is defined by
\[\kt^0\lg:= 0, \qquad \kt^{i+1} \lg := \{ x\in \lg \mid [x,\lg ]\subset  \kt^{i}\lg\text{ and } [Jx,\lg ]\subset  \kt^{i}\lg\}\]
Every subspace $\kt^i\lg$ is $J$-invariant and an ideal of $\lg$. If there is a $k$ such that $\kt^k\lg= \lg$ then the complex structure $J$ is called \emph{nilpotent}.
We have $\kt^i\lg\subset \kz^i\lg$ for all $i$ and if $\kt^i\lg=\kz^i\lg$ then $\kt^{i+1}\lg$ is the largest complex subspace of $\lg$ contained in  $\kz^{i+1}\lg$ (\cite{cfgu00}, Lem. 3).

Not every complex structure is nilpotent, indeed there are examples such that $ \kt^{i}\lg=0$ for all $i$. In general, if $\lg$ is $s$-step nilpotent and the complex structure is nilpotent then we have $\kt^k\lg= \lg$ for some $k\geq s$ and strict inequality is possible.

Cordero, Fernandez, Gray and Ugarte showed the following \cite{cfgu00}:
\begin{prop}
Let $(\lg, J)$ be a nilpotent Lie-algebra with complex structure. Then $J$ is nilpotent if and only if the associated simply connected Lie-group $G$ has the structure of a $k$-step iterated principal $\IC^{n_i}$-bundle where $k$ is the smallest integer such that $\kt^k\lg:= \lg$. The structure of iterated bundle is induced by the filtration $(\kt^i\lg)$ on the Lie-algebra. 
\end{prop}
This bundle structure does not induce a bundle structure on compact quotients of $G$ in general.

%
%

\subsubsection{The complex geometry of $M=\Gamma\backslash G$} We will now address the question of the compatibility of the lattice $\Gamma\subset G$ with the other two structures $\lg$ and $J$. Most of the cited results originate from the work of Malcev \cite{malcev51}.

\begin{defin}
Let $\lg$ be a nilpotent Lie-algebra. A \defobj{rational structure} for $\lg$ is a subalgebra $\lg_\IQ$ defined over $\IQ$ such that $\lg_\IQ\tensor \IR =\lg$.

A subalgebra $\lh\subset \lg$ is said to be rational with respect to a given rational structure $\lg_\IQ$ if $\lh_\IQ:=\lh\cap \lg_\IQ$ is a rational structure for $\lh$.

If $\Gamma$ is a lattice in the corresponding simply connected Lie-group $G$ then its associated rational structure is given by the $\IQ$-span of $\log \Gamma$. A rational subspace with respect to this structure is called \defobj{$\Gamma$-rational}.

A complex structure $J$ on $\lg$ is called {rational} if it maps $\lg_\IQ$ to itself; if $\lg_\IQ$ is induced by a lattice $\Gamma$ in the corresponding Lie-group  $J$ is called $\Gamma$-rational. 
\end{defin}

\begin{rem}\label{ratstrrem}
The lattice and the rational structure are closely related:
\begin{enumerate}
\item The rational structure associated to $\Gamma$ is in fact a rational structure (\cite{Cor-Green}, p. 204) and there exists a lattice in a nilpotent simply connected Lie-group $G$ if and only if the corresponding Lie-algebra admits a rational structure.
\item After passing to a subgroup of finite index one can always assume that  $\log(\Gamma)\subset \lg$ is a lattice, i.e., an additive subgroup which is closed under the Lie bracket. Geometrically this corresponds to a finite covering.
\item If two lattices $\Gamma_1, \Gamma_2\subset G$ generate the same rational structure then the lattice $\Gamma_1\cap\Gamma_2$ is of finite index in both of them (\cite{Cor-Green}, Theorem 5.1.12, p. 205). Geometrically this corresponds to an \'etale correspondence
\[\xymatrix{ &(\Gamma_1\cap\Gamma_2)\backslash G\ar[dl]\ar[dr]\\
\Gamma_1\backslash G && \Gamma_2\backslash G.}\]
\item The subalgebras $\kc^{i}\lg$ and  $\kz^{i}\lg$ as defined above are  always rational subalgebras (\cite{Cor-Green}, p. 208). Thus the filtration $(\kz^{i}\lg)$ yields a structure of real analytic iterated principal torus bundle on $\Gamma \backslash G$.
\end{enumerate}
\end{rem}

\begin{lem}
If in the situation of Lemma \ref{fibration} $\Gamma \subset G$ is a lattice then $\pi$ induces a fibration on the compact nilmanifold $M=\Gamma\backslash G$ if and only if $\Gamma\cap H$ is a lattice in $H$  if and only if the associated subalgebra $\lh\subset\lg$ is $\Gamma$-rational.

In particular $M$ fibres as a principal holomorphic torus bundle $\pi:M\to M'$ over some nilmanifold $M'$ if and only if there is a $J$-invariant, $\Gamma$-rational subspace contained in the centre $\kz\lg$ of $\lg$.
\end{lem}
\pf The first part is \cite{Cor-Green}, Theorem. 5.1.11, p. 204 and Lem 5.1.4, p. 196. The second assertion then follows. \qed

To describe effectively the how the geometry of the nilmanifold is encoded in the Lie-algebra we need
\begin{defin}\label{stableseries}
Let $\lg$ be a nilpotent Lie-algebra with rational structure $\lg_\IQ$. We call an ascending filtration $(\ks^i\lg)_{i=0,\dots, t}$ on $\lg$ a \defobj{torus bundle series} for a complex structure $J$ if for all $i=1\dots , t$
\begin{gather*}
\ks^i\lg \text{ is rational with respect to $\lg_\IQ$ and an ideal in }\ks^{i+1}\lg, \tag{$a$}\\
J\ks^i\lg=\ks^i\lg,\tag{$b$}\\
\ks^{i+1}\lg/\ks^{i}\lg \text{ is an abelian  }.\tag{$c$}
\intertext{If in addition}
\ks^{i+1}\lg/\ks^{i}\lg\subset\kz(\lg/\ks^{i}\lg),\tag{$c'$}\label{princ}
\end{gather*}
then $(\ks^i\lg)_{i=0,\dots, t}$ is called a \defobj{principal torus bundle series}.

An ascending filtration $(\ks^i\lg)_{i=0,\dots, t}$ on $\lg$ is said to be a  \defobj{stable torus bundle series (STBS)}  for $\lg$, if $(\ks^i\lg)_{i=0,\dots, t}$  is a  torus bundle series for every complex structure $J$ and every rational structure $\lg_\IQ$ in $\lg$. If also condition \refb{princ} holds then it is called a  \defobj{stable principal torus bundle series (SPTBS)}.
\end{defin}
Note that we always refer to complex torus bundles and not to real torus bundles. The rationality condition holds for example if every subspace in the series can be described by the subspaces of the ascending and descending central series. 

\begin{rem} If $M=(\lg, J, \Gamma)$ is a nilmanifold and $\lg$ admits a torus bundle series $(\ks^i\lg)_{i=0,\dots, t}$ compatible with $J$ and $\Gamma$ then $M$ can be viewed as an iterated bundle in the following sense: every successive quotient $0\to \ks^{i}\lg\to \ks^{i+1}\lg\to \ks^{i+1}\lg/\ks^{i}\lg\to 0$ induces a holomorphic fibration over a torus with fibre a nilmanifold with left-invariant complex structure. There is not necessarily a iterated  torus bundle structure as discussed before since we did not restrict to the case where $\ks^{i+1}\lg/\ks^{i}\lg$ is an ideal in $\lg/\ks^{i}\lg$; the map $\lg/\ks^{i}\lg \to \lg/\ks^{i+1}\lg$ might not induce a fibration on the nilmanifold. If $J$ is $\Gamma$-rational then the adapted descending central series $\kc^i_J\lg$ gives such an example.

The most intersting case is when  $(\ks^i\lg)_{i=0,\dots, t}$ is a SPTBS for $\lg$. Condition \refb{princ} implies that   $\ks^{i+1}\lg/\ks^{i}\lg$ is an ideal in $\lg/\ks^{i}\lg$  so that  the filtration $(\ks^i\lg)_{i=0,\dots, t}$ induces a  real analytic iterated principal torus bundle structure on $M$:
\[\xymatrix{ T_1 \ar@{^(->}[r] & M_1=M \ar[d]^{\pi_1}\\
T_2 \ar@{^(->}[r] & M_2\ar[d]^{\pi_2}\\
&\vdots\ar[d]\\
T_{s-1}\ar@{^(->}[r] & M_{t-1}\ar[d]^{\pi_{t-1}}\\
& M_t=T_t.}\]
By definition the maps $\pi_i$ are holomorphic with respect to any left-invariant complex structure $J$ because the differentials at the identity are given by the complex linear maps
\[\ks^{i}\lg/\ks^{i-1}\to \lg/\ks^{i-1}.\]
A SPTBS gives in some sense the appropriate description of the complex geometry of the corresponding nilmanifolds.

In general, a STBS, let alone a SPTBS,  will not exist for a nilpotent Lie-algebra but it is not as uncommon as one might think at first sight. In real dimension 6 the Lie-algebra $\lh_7$ described in Example \ref{badex} is the only one not to admit at least a STBS.
\end{rem}
%
\subsubsection{Abelian and complex parallelisable structures}
There are two extremal types of complex structures that can occur in a nilpotent Lie-algebra with complex structure $(\lg, J)$. If $[\einsnull \lg, \nulleins\lg ]=0$ then we have $[Jx.y]=J[x,y]$ for all $x, y\in \lg$; thus the Lie-bracket is $\IC$-linear and $(\lg,J)$ is a complex Lie-algebra.
Nilmanifolds of type $(\lg, J)$ are then complex parallelisable. These have been studied in detail by Winkelmann \cite{winkelmann98} and have very special (arithmetic) properties. Both the descending and the ascending central series are principal torus bundle series (which do not need to be stable).
Their small deformations have been described in \cite{rollenske08a}. 

If $(\lg, J)$ satisfies the complementary condition $[\einsnull{\lg}, \einsnull\lg]=0$ then the complex structure is called abelian and for such nilmanifolds the ascending central series is a  principal torus bundle series.
These structure were extensively studied by Salamon, McLaughlin, Pedersen, Poon, Console, Fino e.a. for example in \cite{con-fin-poon06} and \cite{mpps06}.

\subsection{Dolbeault cohomology of nilmanifolds and small deformations}\label{smalldefo}

Let $M_J=(\lg,J, \Gamma)$ be a nilmanifold with left-invariant complex structure. In this section we want to describe how cohomology and deformations of $M_J$ are governed by $(\lg, J)$.

By a theorem of Nomizu \cite{nomizu54} the de Rham cohomology of  a nilmanifold can be calculated using left-invariant differential forms and is isomorphic to the cohomology of the complex 
\[0\to \lg^*\overset{d}{\to}\Lambda^2\lg^* \overset{d}{\to}\Lambda^3\lg^* \overset{d}{\to}\dots,\]
which also calculates the Lie-algebra cohomology of $\lg$.

The question if the Dolbeault cohomology of nilmanifolds with left-invariant complex structure can be calculated using invariant differential forms has been addressed by  Console and Fino in \cite{con-fin01} and  Cordero, Fernandez, Gray and Ugarte in \cite{cfgu00}: consider in the Dolbeault-complex for $\Omega_{M_J}^p$ the subcomplex of left-invariant forms
\[ 0\to \Lambda^p\einsnull{\lg^*}\overset{\delbar}{\to}\Lambda^p\einsnull{\lg^*} \tensor \nulleins{\lg^*} \overset{\delbar}{\to}\Lambda^p\einsnull{\lg^*} \tensor \Lambda^2 \nulleins{\lg^*}\overset{\delbar}{\to}\dots\]
and denote by $H^{p,q}(\lg,J)$ its $q$-th cohomology group.
Using our definitions their results read as follows:
\begin{theo}\label{citedolbeault}
Let $M_J=(\lg, J,\Gamma)$ be a nilmanifold with left-invariant complex structure. Then the inclusion 
\begin{equation}\label{cohomiso}
\iota_J:H^{p,q}(\lg,J)\to H^{p,q}(M_J) 
\end{equation}
is an isomorphism if
\begin{enumerate}
\item The complex structure $J$ is bi-invariant, $G$ is a complex Lie-group  and $M_J$ is complex parallelisable \cite{sakane76}.
\item The complex manifold $M_J$ is an iterated principal holomorphic torus bundle \cite{cfgu00}.
\item  The complex structure $J$ is $\Gamma$-rational. (\cite{con-fin01}, Theorem B)
\item The complex structure $J$ is abelian \cite{con-fin01}.
\item $\lg$ admits a torus bundle series for $J$ compatible with the rational structure induced by $\Gamma$.
\end{enumerate}
Moreover there is a dense open subset $U$ of the space $\kc(\lg)$ of all left-invariant complex structures on $M$ such that $\iota$ is an isomorphism for all $J\in U$ (\cite{con-fin01}, Theorem A).
\end{theo}
Statement \refenum{v} follows immediately from the the proof of Theorem B in  \cite{con-fin01}. 
The idea of their proof is the following: as long as $M_J$ can be given a structure of iterated bundle with a good control over the cohomology of the base and of the fibre, one can use the Borel spectral sequence for Dolbeault cohomology in order to get an inductive proof. If the complex structure is $\Gamma$-rational then the descending central series adapted to $J$ induces such an iterated fibration which yields the result on a dense subset of the space of invariant complex structures; Console and Fino then show that the property \emph{"The map $\iota_J$ is an isomorphism." } is stable under small deformations.

It is an open question if $\iota_J$ is  an isomorphism for every nilmanifold with left-invariant complex structure.

Generalising results on deformations of abelian complex structures \cite{mpps06, con-fin-poon06} we proved the following
\begin{theo}[\cite{rollenske08c}, Theorem 2.6]\label{invariantdeformation}
Let $M_J=(\lg, J, \Gamma)$ be a nilmanifold with left-invariant complex structure such that 
\[\iota:H^{1,q}((\lg,J),\IC)\to H^{1,q}(M_J)\]
 is an isomorphism for all $q$. Then all small deformations of the complex structure $J$ are again left-invariant complex structures. More precisely, the Kuranishi family contains only left-invariant complex structures. 
\end{theo}

Together with Theorem \ref{citedolbeault} \refenum{v} this yields
\begin{cor}\label{stableinvariantdeformation}
Let $M_J$ be a nilmanifold of type $(\lg,\Gamma)$. If $\lg $ admits  a stable (principal) torus bundle series then every small deformation of $M_J$ is a  nilmanifold with left-invariant complex structure of type $(\lg, \Gamma)$.
\end{cor}

In other words, if $M$ is a real nilmanifold with Lie-algebra $\lg$ and $\lg$ admits a STBS then  the subset of left-invariant complex structures is open in the space of all complex structures. This is a first step on on our way to prove that in some cases left-invariant complex structures form a connected component of the space of complex structures.

\subsection{Examples and Counterexamples}\label{examsection}

The simplest example of a nilmanifolds with left-invariant complex structure is  a complex torus which corresponds to an abelian Lie-algebra.
Other well known examples are Kodaira surfaces, which were historically the first manifolds  known to admit both a symplectic structure and a complex structure but no K\"ahler structure. These are principal bundles of elliptic curves over elliptic curves. 

The following list gives an overview over some phenomena that can occur:
\begin{itemize}
\item Small deformations of iterated principal bundles do not need to have any torus bundle structure: Example \ref{badex}.
\item First examples with SPTBS: Complex tori, Kodaira surfaces, principal bundles over elliptic curve (Example \ref{torusovercurvedescription}).
\item The minimal principal torus bundle series must not to be stable: Example \ref{badmtbs}.
\item A Lie-algebra which admits only non-nilpotent complex structures: $\lh_{26}^+$ discussed in Section \ref{dim3}.
\item A Lie-algebra which admits both nilpotent and non-nilpotent complex structures: Example \ref{nilnonnilexam}, see also \cite{con-fin-poon06}, Section 6.
\item The dimension of the Albanese variety is not constant under small deformations: Example \ref{torusovercurvealbanese}.
\end{itemize}

In order to describe the examples we have to explain some notation. Consider a Lie-algebra $\lg$ spanned by a basis $ e_1, \dots , e_n$. Then the Lie-bracket is uniquely determined by structure constants $a_{ij}^k$ such that 
\[ [e_i, e_j]=\sum_{k=1}^n a_{ij}^ke_k\] 
satisfying  $a_{ij}^k=-a_{ji}^k$ and  the relations encoding the Jacobi identity.
Let $\langle e^1, \dots, e^n\rangle$ be the dual basis, i.e. $e^i(e_j)=\delta_{ij}$. Then for any $\alpha\in \lg^*$ and $x,y\in \lg$ we define
\[d:\lg^*\to \Lambda^2\lg^*, \qquad d\alpha(x,y):=-\alpha([x,y])\]
and get a dual description of the Lie-bracket by
\[de^k=-\sum_{i<j}a_{ij}^k e^{ij} \]
where we abbreviate $e^i\wedge e^j=e^{ij}$.
The map $d$ induces a map on the exterior algebra $\Lambda^*\lg^*$ and $d^2=0$ is equivalent to the Jacobi identity.
Giving a 2-step Lie-algebra is the  same as giving a surjective alternating bilinear form on $\lg/\kz\lg$ with values in $\kz\lg$; the Jacobi Identity is trivial in this case.

If we write $\lh_2 = (0,0,0,0,12,34)$  we mean the following:
Let $e_1, \dots , e_6$ be a basis for the Lie-algebra and $e^1, \dots, e^6$ be the dual basis. Then the defining relations for $\lh_2$ are given by
\[de^1=de^2=de^3=de^4=0,\quad de^5=e^{12}:=e^1\wedge e^2, \quad de^6=e^{34}.\]
Dually, the bracket relations are generated by $[e_1,e_2]=-e_5$ and $[e_3, e_4]=-e_6$. This Lie-algebra will come up again in Section \ref{dim3} when we discuss complex structures on real 6-dimensional Lie-algebras.

We will usually describe a complex structure on a Lie-algebra by specifying an endomorphism $J$ on the basis of $\lg$ such that $J^2=-id_\lg$; if we set $Je_1=e_2$ then $Je_2=J^2e_1=-e_1$ is tacitly understood. 
An equivalent way is to give the decomposition $\lg_\IC=\einsnull\lg \oplus \nulleins\lg$ in the $(\pm i)$-eigenspaces of $J$.

There is a very simple class of examples to which our results will apply:
\begin{exam}\label{torusovercurvedescription}
Let  $M_J$ be a  non-trivial principal holomorphic torus bundle over an elliptic curve. Regarding $M_J$ as a nilmanifold with left-invariant complex structure with Lie algebra $\lg$, the fibration corresponds to a central extension
\[0\to \kt^1\lg \to \lg\to \mathfrak{e}\to 0\]
where $\kt^1\lg$ is the $(2n-2)$-dimensional real, $J$-invariant subspace corresponding to the fibres of the bundle.
In particular the centre of $\lg$ has dimension at least $2n-2$ whence, since $M$ is not a torus, we have  $\kz\lg=\kt^1\lg$. There is only one Lie-algebra with 2-codimensional centre so  $\lg\isom (0,0,0,\dots, 0, 12)$. We will see in Proposition \ref{propdimcg=1} that
\[\lg\supset \kz\lg\supset 0\]
 is a SPTBS for $\lg$.
\end{exam}

%

\begin{exam}\label{badex}
A nice behaviour of the complex structure on the universal covering $G$ may not be sufficient to get a good description the manifold $M=\Gamma\setminus G$.

Consider the Lie-group $H_7$ whose Lie-algebra is $\lh_7= (0,0,0,12,13,23)$ with basis $e_1, \dots, e_6$.  The vectors $e_4\dots,e_6$ span the centre $\kz^1\lh_7=\kc^1\lh_7$. 

Let $\Gamma\subset H_7$ be the lattice generated the elements $\exp (e_k)$ and consider the nilmanifold $M=\Gamma\backslash H_7$, which -- as a real manifold -- can be regarded as a real principal torus bundle with fibre and base a 3-dimensional torus.

For every real number $\lambda$  we give a left-invariant complex structure $J_\lambda$ on $M$ by specifying a basis for the space of $(1,0)$-vectors:
\[(\einsnull{\lh_7})_\lambda:=\langle X_1=e_1-ie_2, X_2^\lambda= e_3-i (e_4-\lambda e_1), X_3^\lambda=-e_5-\lambda e_4+ie_6\rangle\]
One can check that $[X_1, X_2^\lambda]=X_3^\lambda$ and, since $X_3^\lambda$ is contained in the centre, the complex structure is integrable. The largest complex subspace of the centre is spanned by the real and imaginary part of $X_3^\lambda $ since the centre has real dimension three.

The simply connected Lie-group $H_7$ has  a filtration by subgroups induced by the filtration 
\[\lh_7\supset V_1=\langle \lambda e_2+e_3, e_4,Im(X_3^\lambda),Re(X_3^\lambda)\rangle\supset V_2=\langle Im(X_3^\lambda),Re(X_3^\lambda)\rangle \supset 0\]
on the Lie-algebra and, since all these are $J$ invariant, $H_7$ has the structure of a tower of principal holomorphic bundles with fibre $\IC$, i.e.,  the complex structure is nilpotent.
In fact, using the results of \cite{math.DG/0411254} a simple calculation shows that every complex structure on $\lh_7$ is equivalent to $J_0$.

Now we take the compatibility with  the lattice into account. 
The rational structure induced by $\Gamma$ coincides with the $\IQ$-algebra generated by the basis vectors $e_k$. Therefore $V_j\cap\Gamma$ is a lattice in $V_j$ if and only if $\lambda$ is rational. That is, for $\lambda\notin \IQ$ the structure of iterated holomorphic principal bundle on $H_7$ does not descend to the quotient $\Gamma\backslash H_7$.
\end{exam}

\begin{exam}\label{badmtbs} The aim of this example is to show that  the minimal torus bundle series $\kt^i\lg$ might not yield the appropriate geometric description of an associated nilmanifold.

We consider the Lie-algebra given by $\lg=(0,0,0,0,0,0,12,13)$ with the lattice $\Gamma\subset G$ generated by the images of the  basis vectors $e_k$.  Let $\langle z_1, z_2, z_3\rangle$ be any basis of the subspace spanned by $e_4, \dots ,e_6$
We define a complex structure by declaring 
\begin{align*}
Je_2&=e_3,& Je_7&=e_8,\\ Jz_1&=e_1,&Jz_2&=z_3.
\end{align*}
A short calculation shows that $J$ is integrable and that $\kt^1\lg=\langle z_2, z_3, e_7, e_8\rangle$ is the largest complex subspace of the centre. But the $z_i$ were quite arbitrarily chosen and hence $\kt^1\lg$ is not a $\Gamma$-rational subspace in general, which means that it does not tell us a lot about the geometry of the compact manifolds $M=\Gamma\setminus G$. But on the other hand $\kc^1\lg=\langle e_7, e_8\rangle$ is a $\Gamma$-rational, $J$-invariant subspace. This is no coincidence: we will  prove in Proposition \ref{dimcg=2} that
\[0\subset \kc^1\lg \subset \lg\]
is a SPTBS for $\lg$ which means that we can describe $M$ as a principal bundle of elliptic curves over a complex 3-torus for any complex structure on $\lg$.
\end{exam}

\section{Albanese-Quotients and deformations in the large}\label{albanesequotients}

In this section $\Delta$ will denote a small 1-dimensional disc centred in 0.
In  \cite{catanese04} Catanese proved that any deformation in the large of a complex torus is again a complex torus by analysing the Albanese map:
he studies families $\kx\to \Delta$ with general fibre a complex torus. In a first step he shows that also the special fibre $\kx_0$ has a surjective Albanese map to a complex torus of the same dimension. By proving that this map has to be biholomorphic he concludes that  the special fibre is a complex torus. 
We will try to generalise his method to the case where the Albanese map is not so well behaved.

\subsection{Definitions and results}

We need to recall some definitions.
Let $X$ be a compact, complex manifold. By Kodaira's Lemma (\cite{catanese04}, Lemma 2.2) we have an inclusion 
\[H^0(X,d\ko_X)\oplus \overline{H^0(X,d\ko_X)}\subset H^1_{dR}(X,\IC).\]

\begin{defin}
The \defobj{Albanese variety} of $X$ is the  abelian, complex  Lie-group $Alb(X)$  defined as the quotient of
$H^0(X,d\ko_X)^*$ by the minimal, closed, complex subgroup containing the image
of $H_1(X, \IZ)$ under the map
\[H_1(X, \IZ)\to H_1(X, \IC)\to H^0(X,d\ko_X)^*.\]
 The Albanese variety is called very good if the image of
$H_1(X, \IZ)$ is a lattice  in $H^0(X,d\ko_X)^*$.
\end{defin}
The Albanese map $\alpha_X:X\to Alb (X)$ is given by integration of closed forms along paths starting from a fixed base point. 

\begin{rem}
J\"org Winkelmann brought to our attention that the Albanese variety is always compact and hence $Alb(X)$ is very good if and only if the image of $H_1(X,\IZ)$ is discrete in $H^0(X,d\ko_X)^*$.

If $X$ satisfies the weak 1-Hodge property which means
\[H^1_{dR}(X,\IC)=H^0(X,d\ko_X)\oplus \overline{H^0(X,d\ko_X)}= H^0(X,d\ko_X)\oplus H^1(X, \ko_X)\]
then $X$ has a very good Albanese variety; in particular this holds if $X$ is k\"ahlerian.
\end{rem}

Catanese studied the behaviour of the Albanese map with respect to deformation and obtained for example the following
\begin{prop}[\cite{catanese04}, Corollary 2.5] Assume that $\{X_t\}_{t\in\Delta}$ is a 1-parameter family of compact complex manifolds over the unit disc such that there is a sequence $t_\nu\to 0$ with $\kx_{t_\nu}$ satisfying the weak 1-Hodge property, and moreover such that the image of $X_{t_\nu}$ under the Albanese map has dimension $a$.
Then also the central fibre has a very good Albanese map and Albanese dimension $a$.
\end{prop}

In the following we want to explain how to generalise this result if the Albanese variety is not necessarily very good; the right way to do so is best explained via the following example.

\begin{exam}\label{torusovercurvealbanese}
Let $X=(\lg, J, \Gamma)$ be a  principal holomorphic torus bundle over an elliptic curve  as in  Example \ref{torusovercurvedescription}. Then there is a  basis $e_1, e_2, z_1, \dots , z_{2k}$ of $\lg$ such that the only non-trivial Lie-bracket relation is $[e_1, e_2]=z_{2k}$ and $\kc^1\lg$ is spanned by $z_{2k}$. The centre $\kz\lg=\langle z_1, \dots, z_{2k}$ corresponds to the fibres of the bundle and is always $J$-invariant.
With respect to the dual basis $e^1, \dots, z^{2k}$  the differential $d:\lg^*\to \Lambda^2\lg^*$ is given  by 
\[ dz^{2k}=-e^1\wedge e^2, de^1=de^2=dz^1=\dots=dz^{2k-1}=0\]
and, by the description of cohomology given in Section \ref{smalldefo}, we see that
\begin{gather*}
H^1(X,\IC)=\{\alpha \in \lg^*_\IC\mid d\alpha =0\}
\end{gather*}
has dimension $2k+1$. Thus $X$ can never satisfy the weak 1-Hodge property.

In order to determine the Albanese Variety of $X$ we have to take the rational structure on $\lg$ induced by the lattice $\Gamma$ into account;  on the quotient $H_1(X,\IR)=\lg/\kc^1\lg$ it maps to  the rational structure induced by $H_1(X,\IZ)$.
Let $W$ be the smallest $J$-invariant, $\Gamma$-rational subspace of $\lg$ which contains $\kc^1\lg$. Since the centre is $J$-invariant and $\Gamma$-rational we have
\[2\leq \dim_\IR W \leq \dim_\IR \kz\lg=2k.\]

The space of closed holomorphic 1-forms is
\begin{align*} H^0(X,d\ko_X)&=\{\alpha \in \einsnull{\lg^*}\mid d\alpha=0\}\\
&= \{\alpha \in \einsnull{\lg^*}\mid\delbar\alpha=0\}\\
&= Ann((\kc_J^1\lg)\tensor\IC)\cap \einsnull{\lg^*},
\end{align*}
 hence the kernel of the surjection $p:H_1(X,\IR)\onto H^0(X, d\ko_X)^*$ is  $\kc_J^1\lg/\kc^1\lg$ yielding
\[Alb(X)=\frac{H^0(X, d\ko_X)^*/p(W)}{\im (H_1(X,\IZ)\to H^0(X, d\ko_X)^*/p(W))}. \]

The dimension of the Albanese variety can vary with the complex structure in the range
\[1\leq\dim Alb(X)=\frac{\dim_\IR \lg-\dim_\IR W}{2}\leq k= h^0(X,d\ko_X).\]

But in some sense all the ugly phenomena happen inside $\kz\lg$. Indeed, if we replace $W$ by $\kz\lg$ in the above construction we get a map to an elliptic curve such that the corresponding map in cohomology is independent of the complex structure; we have found a way to reconstruct our map $\pi:X\to E$ from cohomological data which does not depend on the particular left-invariant complex structure.
\end{exam}

Thus, even if the Albanese map is not well behaved, we might find some subspace in $H_1(X,\IR)$, not depending on the complex structure, which yields a holomorphic map to a compact torus. We will now formalise this idea.

Let $V\subsetneq H_1(X,\IQ)$ be a proper subspace and let 
\[U:=Ann(V_\IC)\cap H^0(X,d\ko_X)\]
where the annullator is taken with respect to the pairing between $H^1(X, \IC)$ and $H_1(X,\IC)$.
If $Ann(V_\IC)=U\oplus \bar U$ then the map $H_1(X, \IR)/V_\IR\overset{\isom}{\to} U^*$ is an isomorphism of real vector spaces.

Since $V$ is a rational subspace, $H_1(X, \IZ)/(H_1(X, \IZ)\cap V)$ maps to a full lattice  $\Lambda\subset U^*$.
This yields a complex torus $U^*/\Lambda$ which is a quotient of the Albanese variety $Alb(X)$, the projection being induced by the inclusion  $U\into  H^0(X,d\ko_X)$.

We have shown that the following is well defined and yields in fact a compact, complex torus:
\begin{defin}\label{qalbdefin} Let $V\subsetneq H_1(X,\IQ)$ be a proper subspace such that
\[Ann(V_\IC)=U\oplus \bar U\]
where $U:=Ann(V_\IC)\cap H^0(X,d\ko_X)$.
Then we call 
\[ QAlb_V(X):=U^*/\Lambda\]
the \defobj{very good Albanese-Quotient} of $X$ associated to $V$, where $\Lambda$ is the image of the lattice $H_1(X, \IZ)/(H_1(X, \IZ)\cap V)$
under the map $H_1(X, \IR)/V_\IR\to U^*$.
\end{defin}

\begin{rem}\label{Vcond}
The above condition on $V$ is satisfied if and only if the kernel of the composition map
\[\phi: H^0(X,d\ko_X)\to H^1(X, \IC)\to H^1(X, \IC)/Ann(V_\IC)\]
has complex dimension $q=\frac{\dim_\IQ(V)}{2}$. This is in fact the maximal possible dimension since $V_\IC$ is defined over $\IR$ while on the other hand
\[H^0(X,d\ko_X)\cap \overline{H^0(X,d\ko_X)}=0.\]

Another way to look at this is the following: let
\[W:= (H^0(X,d\ko_X)\oplus \overline{H^0(X,d\ko_X)})\cap H^1(X, \IR).
\]
The decomposition $W_\IC =H^0(X,d\ko_X)\oplus \overline{H^0(X,d\ko_X)}$ defines a complex structure on $W$ and the above conditions on $V$ are equivalent  to $Ann(V_\IR)$ being a complex subspace of $W$ with respect to this complex structure.
\end{rem}

\begin{rem}
Integration over the closed forms which are contained in $Ann(V_\IC)\cap H^0(X,d\ko_X)$ gives us a holomorphic map $q\alpha_V$ which factors through the Albanese map:
\[\xymatrix{X\ar[r]^\alpha \ar[dr]_{q\alpha_V} & Alb (X) \ar[d]\\ & QAlb_V(X).}\]

Note that $QAlb_V(X)$ has positive dimension by definition since we assumed $V\neq H_1(X,\IQ)$.
\end{rem}

\begin{exam}
Let $M=(\lg, J, \Gamma)$ be a 2-step nilmanifold such that the centre $\kz\lg$ is $J$-invariant. By Nomizu's theorem we have $H_1(M,\IR)\isom \lg/[\lg, \lg]$ and $V_\IR:=\kz\lg/[\lg, \lg]$ is a $\Gamma$-rational subspace. Then $V_\IQ$ induces a very good Albanese-Quotient and $q\alpha_V$ coincides with the principal holomorphic bundle map induced by the filtration $\lg\supset\kz\lg\supset 0$. In particular this applies to the principal holomorphic torus bundle over an elliptic curve discussed in Example \ref{torusovercurvealbanese} above.
\end{exam}

Our goal is to study such locally trivial, smooth, holomorphic fibrations. Hence we define
\begin{defin}\label{goodfibreclass}
Let $\kc$ be a class of compact, complex manifolds and let $X$ be any compact, complex manifold. We say that a subspace $V\subsetneq H_1(X,\IQ)$ is \defobj{$\kc$-fibering} on $X$ if $V$ induces a very good Albanese-Quotient and the map
\[q\alpha:X\to QAlb_V(X)\]
is a locally trivial, smooth, holomorphic fibration with fibres  in $\kc$.

A class $\kc$ of compact complex manifolds is called a \defobj{good fibre class} if the following conditions hold
\begin{enumerate}
\item The class  $\kc$ is closed under holomorphic 1-parameter  limits, i.e., if we have a smooth family over the unit disc  and a sequence $(t_\nu)_{\nu\in \IN}$ converging to 0 such that the fibres over $t_\nu$ are in $\kc$ then also the central fibre is in $\kc$.
\item One of the following conditions holds:
\begin{enumerate}
\item There is a coarse moduli space $\mathfrak M_\kc$ for manifolds in class $\kc$ which is Hausdorff.
\item For every manifold  $Y\in\kc$ there is a local moduli space, i.e.,  the Kuranishi family of $Y$ is universal.
\item $h^1(F, \Theta_F)$ is constant on the connected components of $\kc$.
\end{enumerate}
\end{enumerate}
\end{defin}

So in the above example the subspace $V$ is torus-fibering. 

\begin{rem}\label{C-rem}
The condition (a) is very strong but it is worth including since it applies to Kodaira surfaces. Indeed these and complex tori are our first examples of good fibre classes.

That every deformation in the large of a complex torus is again a complex torus has been proved by Catanese in \cite{catanese02} and we have $h^1(T, \Theta_T)=\dim(T)^2$ for any Torus $T$ (for example by Section \ref{smalldefo}).

Borcea proved in  \cite{borcea84} that the moduli space of Kodaira surfaces may be identified with the product of the complex plane by the punctured disc which is Hausdorff. That every deformation in the large is again a Kodaira surface has been proved by the author in  \cite{rollenske05b}.
\end{rem}

The main result of this section is

\begin{theo}\label{largedeform}
Let $\kc$ be a good fibre class and $\pi:\kx\to \Delta$ be a smooth family of compact, complex manifolds. Let $(t_\nu)_{\nu\in \IN}$ be a sequence in $\Delta$ converging to zero.

If there is a subspace $V\subsetneq H_1(\kx_0,\IQ)$ which is $\kc$-fibering on $\kx_{t_\nu}$ for all $\nu$ then $V$ is $\kc$-fibering on $\kx_0$. 
\end{theo}

The proof will be given in the next subsection. We will first show how this can be applied to nilmanifolds

\begin{theo}\label{nillargedeform}
Let $G$ be a simply connected nilpotent Lie-group with   Lie-algebra $\lg$ and let $\Gamma\subset G$ be a lattice such that the following holds:
\begin{enumerate}
\item $\lg$ admits a STBS  $(\ks^i\lg)_{i=0,\dots, t}$ (cf. Definition \ref{stableseries}).
\item The nilmanifolds with left-invariant complex structure of the type $(\ks^{t-1}\lg, \Gamma\cap \exp(\ks^{t-1}\lg))$ constitute a good fibre class.
\end{enumerate}
Then any deformation in the large of a nilmanifold with left-invariant complex structure  $M_J=(\lg,J,\Gamma)$ is again of type $(\lg, \Gamma)$.
\end{theo}

\pf We will use the following:
\begin{lem}[\cite{catanese04}, Lemma 2.8] Let $\kb$ be a connected complex analytic space and $\kb'$ a non-empty, open subset of $\kb$ such that $\kb'$ is closed for holomorphic 1-parameter limits (i.e., given any holomorphic map $f:\Delta \to \kb$ , if there is a sequence $(t_\nu)_{\nu\in \IN}$ converging to 0 with $f(t_\nu)\in \kb'$ then $f(0)\in \kb'$). Then $\kb=\kb'$.
\end{lem}

Now consider any smooth family $\kx\to \kb$ over a connected base $\kb$ such that for some point $b\in \kb$ the fibre $\kx_b$ is isomorphic to a nilmanifold $M=(\lg, J, \Gamma)$. Then consider the (non-empty) set
\[ \kb':=\{t\in \kb \mid \kx_t\isom M'=(\lg, J', \Gamma) \text{ for some complex structure }J'\}.\]
This set is open in view of our result on small deformations in Corollary \ref{stableinvariantdeformation}.
We have to show that it is closed under holomorphic 1-parameter limits.

Our definition of STBS guarantees that $\Gamma':=\Gamma\cap \exp(\ks^{t-1}\lg)$ is a lattice. Taking 
\[V_\IR=\ks^{t-1}\lg/\kc^1\lg\subset \lg/\kc^1\lg\isom H^1(M, \IR)\]
we can apply Theorem \ref{largedeform} which yields that any 1-parameter limit $\hat M$ of manifolds in $\kb'$ is a smooth, locally trivial, holomorphic fibration over a torus with fibre a nilmanifold $F=(\ks^{t-1}\lg, J, \Gamma')$. The topological structure of the fibre bundle is determined by the fundamental group, which does not change under deformation, and the fundamental group $\Gamma$ determines $G$ and $\lg$ up to isomorphism (\cite{VinGorbShvart}, p. 45, corollary 2.8). So  also $\hat M$ is a nilmanifolds with left-invariant complex structure of the same type.

We have proved that $\hat M \in \kb'$ hence  $\kb'$ is closed under holomorphic 1-parameter limits and, by the lemma, $\kb=\kb'$ which concludes the proof.\qed \margincom{give example-application?}

Using Remark \ref{C-rem} we see as a first application that every deformation in the large of a principal holomorphic torus bundle over an elliptic curve is again such a principal holomorphic bundle. Further applications will be given in Section \ref{applications}.

\subsection{Proof of Theorem \ref{largedeform}}

We will split the proof of the theorem into several steps.

Catanese showed in \cite{catanese91} that the Albanese dimension is in fact a topological property if $X$ is a K\"ahler manifold and we review his arguments in our context:

Assume that we have $X$ and $V\subset H_1(X,\IQ)$ as in Definition \ref{qalbdefin}.

\begin{lem}\label{dimensionlemma}
The dimension  of the image of $X$ under the map $q\alpha$ is 
\[d=max \{ m\mid \im(\Lambda^m (Ann(V_\IC)\cap H_0(d\ko_X))\to H^m(X,\IC))\neq 0\}.\]
Moreover $q\alpha$ is surjective if and only if  $\Lambda ^{2k} Ann(V_\IC)\neq 0$ in $H^{2k}(X,\IC)$ where $k=\dim Ann(V_\IC)$. In particular, the surjectivity of $q\alpha$ is a property which depends only on the topology of $X$ and the subspace $V$.
\end{lem}
\pf
By the definition of $QAlb_V(X)$  we have
\[U:=q\alpha^*H^0(QAlb_V(X), d\ko_{QAlb_V(X)})= Ann(V)\cap H_0(d\ko_X)\]
and $q\alpha$ is given by integration over the holomorphic 1-forms in $U$.

The dimension of its image is, by Sard's theorem, equal to the maximal rank of the differential of the quotient Albanese map and hence equal to $d$. 

It remains to show  that $\Lambda ^{2k} Ann(V)\neq 0$ if and only if $\Lambda^k U\neq 0$. But this is clear since our assumptions guarantee that $Ann(V)$ is contained in $H^0(X,d\ko_X)\oplus \overline{H^0(X,d\ko_X)}$, hence $Ann(V)_\IC=U\oplus \bar U$ and $\Lambda ^{2k} Ann(V)=\Lambda^k U\tensor \Lambda^k\bar U$.\qed

In order to analyse how our notion of Albanese-Quotient behaves under deformation we need to introduce some notation.

Let $\pi:\kx\to \Delta$ be a smooth family of compact complex manifolds over the unit disk. We identify  $H_1(\kx_t, \IQ)$ with $H_1(\kx_0, \IQ)$.

 The sheaves of relative differential forms $\Omega^k_{\kx/\Delta}$ are defined by the exact sequence 
\[ 0\to \pi^*\Omega^k_{\Delta}\to \Omega^k_{\kx}\to \Omega^k_{\kx/\Delta}\to 0\]
and we have a $\pi^*\ko_\Delta$-linear map 
\[d_v:\Omega^{k-1}_{\kx/\Delta}\overset{d}{\to} \Omega^k_{\kx}\to \Omega^k_{\kx/\Delta}\]
given by differentiation along the fibres.
This gives us a complex of sheaves
\[\ke^\bullet= 0\to \ko_\kx\to \Omega^1_{\kx/\Delta} \to \Omega^2_{\kx/\Delta}\to \dots\]
which restricts to the holomorphic de Rham complex on every fibre $\kx_t$ of $\pi$.

The following useful fact has been extracted from the proof of Lemma 2.4 in \cite{catanese04}. Denote by $\IC_t$ the residue field at the point $t\in \Delta$.
\begin{lem}\label{limitlemma}
Let $\pi:\kx\to \Delta$ be a smooth family of compact complex manifolds over a small disk and assume that there is a sequence $(t_\nu)_{\nu\in \IN}$ converging to 0 such that $H^0(\kx_{t_\nu}, d\ko_{\kx_{t_\nu}})$ has dimension $q$ for all $\nu$.

Then, after possibly shrinking $\Delta$, there is a rank $q$ locally free subsheaf $\kh$ of $H^1(\kx_0, \IC)\tensor \ko_\Delta$,  the trivial vector bundle on $\Delta$ with fibre $H^1(\kx_0, \IC)$, such that 
\[\kh\tensor \IC_t\subset H^0(\kx_t, d\ko_{\kx_t}) \text{ for all }t\in \Delta\]  and equality holds for $t\neq 0$.
\end{lem}
\pf Let $\Delta^*=\Delta\setminus\{0\}$ be the pointed disk. We define a (not necessarily locally free) sheaf $d\ko_{\kx/\Delta}$ on $\kx$ by the exact sequence
\[ 0\to d\ko_{\kx/\Delta}\to \Omega^1_{\kx/\Delta}\overset{d_v}{\to} \Omega^2_{\kx/\Delta}\]
and push down this sequence to $\Delta$ via $\pi_*$. 
 
Since the sheaves $\Omega^k_{\kx/\Delta}$ are torsion free the same holds for their direct images $\pi_*\Omega^k_{\kx/\Delta}$, and, $\Delta$ being smooth of dimension 1,  the $\pi_*\Omega^k_{\kx/\Delta}$ are in fact locally free. This implies that also $\kh:=\pi_*d\ko_{\kx/\Delta}$ is locally free since it is a subsheaf of a locally free  sheaf on $\Delta$. Any section in $\kh$ is a holomorphic differential form on $\kx$ which restricts to a closed form on any fibre.

The base change map
\[\pi_*\Omega^k_{\kx/\Delta}\tensor \IC_t\into H^0(\kx_t, \Omega^k_{\kx/\Delta}\restr{\kx_t})=H^0(\kx_t, \Omega^k_{\kx_t})\]
is an injection, since $\{t\}$ has codimension 1 in $\Delta$ (\cite{Gr-Rem}, Prop. 2, p. 208, p. 209).
By possibly shrinking our disk we may assume that on $\Delta^*$ the dimensions of $H^0(\kx_t, \Omega^k_{\kx_t})$ are constant ($k=1,2$) and  the map of vector bundles $d_v:\pi_*\Omega^{1}_{\kx/\Delta}\to \pi_*\Omega^2_{\kx/\Delta}$ has constant rank.
Then also $h^0(\kx_t,d\ko_{\kx_t})$ is constant on $\Delta^*$ and we have an isomorphism \[\pi_*d\ko_{\kx/\Delta}\tensor \IC_t\isom H^0(\kx_t, d\ko_{\kx_t})\]
 for $t\neq 0$. In particular $\kh$ has rank $q=h^0(\kx_{t_\nu}, d\ko_{\kx_{t_\nu}})$.

The map $\kh\into H^1(\kx_0, \IC)\tensor \ko_\Delta$ is induced by the inclusion
\[H^0(\kx_t, d\ko_{\kx_t})\into H^1(\kx_0, \IC)\]
 on each fibre.\qed

Considering a product of elliptic curves deforming to a simple torus it is clear that a very good Albanese-Quotient may not be stable under small deformations. On the other hand we have:

\begin{prop}\label{Albanesesubset}
Let $\pi:\kx\to \Delta$ be a smooth family of compact, complex manifolds over the unit disk and let $V\subset H_1(\kx_0, \IQ)$. Then the set
\[ Q:=\{ t\in \Delta\mid \text{$V$ defines a very good Albanese-Quotient on $\kx_t$}\}\]
is a (possibly empty) analytic subset of $\Delta$.
\end{prop}

\pf The question is local on $\Delta$ and hence we may assume that we are in the situation of Lemma \ref{limitlemma}. Consider the composition map of vector bundles on $\Delta$ given by
\[\phi:\kh\to H^1(\kx_0, \IC)\tensor \ko_\Delta\to \left(H^1(\kx_0, \IC)/Ann(V_\IC)\right)\tensor \ko_\Delta.\]
If $\dim_\IC Ann(V_\IC)=2q$ then $V$ induces a very good Quotient Albanese map on $\kx_t$ if and only if the kernel of the map $\phi_t$ has dimension $q$ which is the maximal possible dimension (see Remark \ref{Vcond}). Writing $\phi$ as a matrix with holomorphic entries we see that this is equivalent to the vanishing of the determinants of all minors of a certain dimension which is an analytic condition.\qed

\begin{cor}\label{albaneselimit}
Let $\pi:\kx\to \Delta$ be a smooth family of compact, complex manifolds over the unit disk and let $(t_\nu)_{\nu\in \IN}$ be a sequence  converging to some point in $\Delta$. If $V\subset H_1(\kx_0, \IQ)$ gives rise to a very good Albanese-Quotient on $\kx_t$, for $t={t_\nu}$ then this holds for all  $t\in \Delta$.

If the quotient Albanese map is surjective for some $t_\nu$ then it is surjective for all $t$.
\end{cor}

\pf Consider the set $Q$ as in the Proposition. An analytic subset which has an accumulation point must have dimension at least 1. But since $\Delta$ itself is 1-dimensional we have $Q=\Delta$ and the first claim is proved. The last statement follows immediately from \ref{dimensionlemma}.\qed

We can now study to what extend properties of the general fibre transfer to the central fibre. First we show that there is indeed a family of fibrations.

\begin{lem}\label{quotientfamily}
Let  $\pi:\kx\to \Delta$ be  a smooth family over a small disc and let  $V\subset H_1(\kx_0, \IQ)$ be a subspace inducing a very good Albanese-Quotient on every fibre. Then after possibly shrinking $\Delta$ there is a family of tori $\pi': \kb\to \Delta$ and a map $\Phi$ inducing a diagram
\[\xymatrix{ \kx\ar[rr]^\Phi\ar[dr]_\pi &&\kb\ar[dl]^{\pi'}\\ & \Delta}\]
such that for every $t\in \Delta$ the map $\Phi_t:\kx_t\to \kb_t$ is  the quotient Albanese map. 
\end{lem}
We say that $\Phi$ is a family of Albanese-Quotients.

\pf
We may assume that we are in the situation of Lemma \ref{limitlemma} and that there is a local cross section $s:\Delta \to\kx$. We define $\kh':=Ann(V_\IC)\cap \kh$ and get our family $\kb\to \Delta$ by taking the quotient of ${\kh'}^*$ by the image of $H_1(X, \IZ)/(H_1(X, \IZ)\cap V)$. The map $\Phi$ can be defined by mapping a point $x\in \kx$ to the  map $\omega\mapsto \int_{\gamma_x} \omega$ where $\gamma_x$ is any path joining $x $ to $s(\Delta)$. Then $\Phi$ restricts to the quotient Albanese map on every fibre as claimed.\qed

\begin{prop}\label{smoothfibration}
Let 
\[\xymatrix{ \kx\ar[rr]^\Phi\ar[dr]_\pi &&\kb\ar[dl]^{\pi'}\\ & \Delta}\]
be a family  of Albanese-Quotients such that $\Phi_t:\kx_t\to \kb_t$ is a smooth holomorphic fibration for $t=t_{\nu}$ where $(t_\nu)_{\nu\in \IN}$ is a sequence converging to 0 in $\Delta$. Then there is a small neighbourhood $\Delta'$ of zero such that $\Phi_t$ is a smooth holomorphic fibration for all $t\in \Delta'$.
\end{prop}
\pf 
By Proposition \ref{albaneselimit} the map $\Phi_0$ is also surjective.
Recall that $\Phi$ is given by integration over the closed holomorphic 1-forms in
\[U:=Ann(V_\IC)\cap H_0(d\ko_\kx).\]

Let $\omega_1, \dots , \omega_m$ be a basis for $U$. Then $\omega:=\omega_1\wedge\dots\wedge \omega_m$ generates a sub line bundle of $\Omega^m_\kx$, namely the pullback $\Phi^*K_{\kb/\Delta}$ of the relative canonical bundle of the family $\kb\to \Delta$.

The rank of the Jacobian of $\Phi$ in some point $p\in\kx$ is not maximal if and only if $\omega$ vanishes in $p$, i.e. in the points of the zero divisor $R:=Z(\omega)$. But since $\Phi_t$ is smooth for $t=t_\nu$ the divisor $R$ is completely contained in a union of fibres and we can choose $\Delta'$ such that $\Phi_t$ is smooth for $t\in \Delta'\setminus \{0\}$. 
Hence we may assume that $R$ is completely contained in $\kx_0$ and, since $\kx_0$ is irreducible of codimension 1, we have in fact $R=\kx_0$ or $R=\emptyset$.

If $R$ is not empty then there is a maximal $k\in \IN$ such that $\omega/t^k$ is holomorphic and hence, after a base change $\Delta\overset{t^k}{\to} \Delta$, there is at least one point in $\kx_0$ where $\omega$ does not vanish. But then it can vanish nowhere by dimension reasons  and this proves that the Jacobian has maximal rank everywhere and the central fibre is indeed a smooth holomorphic fibration.\qed

Now we want to know, when fibration on the central fibre is locally free provided this holds for $\Phi_t$ $(t\neq 0)$.
 For a family of manifolds $\Phi:\kx\to \kb$  we define the sheaf of relative tangent vectors $\Theta_{\kx/\kb}$ by the sequence
\[ 0\to \Theta_{\kx/\kb}\to  \Theta_\kx  \to \Phi^*\Theta_\kb\to 0.\]

\begin{prop}\label{locallytrivial}
Let\[\xymatrix{ \kx\ar[rr]^\Phi\ar[dr]_\pi &&\kb\ar[dl]^{\pi'}\\ & \Delta}\] be a family of smooth holomorphic fibrations over tori, parametrised by the unit disc, such that all the fibres of the maps $\Phi_t:\kx_t\to \kb_t$ are contained in a class $\kc$ of compact complex manifolds.

If $\Phi_t$ is locally  trivial for $t\neq 0$ then also $\Phi_0$ is locally trivial if one of the following conditions holds:
\begin{enumerate}
\item There is a coarse moduli space $\mathfrak M_\kc$ for manifolds in class $\kc$ which is Hausdorff.
\item For every manifold in $Y\in \kc$ there is a local moduli space, in other words, the Kuranishi family of $Y$ is universal.
\item The sheaf $R^1\Phi_* \Theta_{\kx/\kb}$ is locally free.
 (This holds if $h^1(\inverse\Phi(p), \Theta_{\kx_p})$ does not depend on $p\in \kb$.)
\end{enumerate}
\end{prop}

\pf
In the first case let $F$, $F'$ be two fibres of $\Phi_0$. By \cite{fi-gr65} it suffices to show that $F\isom F'$. By pullback we can obtain two families $\kf$, $\kf'$ of manifolds in $\kc$ parametrised by $\Delta$ such that $\kf_t$ is a fibre of $\Phi_t$ and with central fibres $F$ and $F'$ respectively. Since $\Phi_t$ is locally trivial for $t\neq 0$ by assumption the corresponding moduli maps to $\mathfrak M_\kc$ coincide for $t\neq 0$ and since $\mathfrak M_\kc$ was assumed to be Hausdorff they have to coincide also for $t=0$ thence $F\isom F'$ as claimed.

For the second case we we consider the family $\kx\to \kb$ in the neighbourhood of some point $p_0\in \kb_0$. 

Let $\mathrm{Kur}(\inverse\Phi(p_0))$ be the Kuranishi family of $\inverse\Phi(p_0)$ which is universal by assumption, i.e., in some neighbourhood $U$ of $p_0$ we get a unique moduli map $\mu:U\to \mathrm{Kur}(\inverse\Phi(p_0))$ such that the point $\mu(p)$ corresponds to the manifold $\inverse\Phi(p)$. 

If we choose $U$ sufficiently small we can find local coordinates $x=(t,y)$ around $p_0$ such that $\pi'(x)=t$. Since $\Phi_t:\kx_t\to \kb_t$ is locally trivial for $t\neq 0$ we have
\[\frac{\del \mu(t,y)}{\del y}=0\]
on the dense set where $t\neq 0$ and hence everywhere.

Therefore the moduli map is constant on $U\cap \kb_0$ and, since the point $p_0\in \kb_0$ was arbitrary and $\kb_0$ is connected, all the fibres of $\phi_0$ are isomorphic and  $\Phi_0$ is locally trivial by \cite{fi-gr65}.

It remains to treat the last case.
The fibration $\Phi_0:\kx_0\to \kb_0$ is locally trivial if and only if the Kodaira-Spencer map
\[\rho_0:\Theta_{\kb_0}\to R^1{\Phi_0}_*\Theta_{\kx_0/\kb_0}\]
vanishes identically. (\cite{Gr-Rem}, Proposition 1, p. 208).

We want to study the relation between $\rho_0$ and the Kodaira-Spencer map $\rho$ for the whole family via the base change homomorphism. 
In particular we are interested in the subsheaf $\Theta_{\kb/\Delta} \subset \Theta_\kb$ that restricts to $\Theta_{\kb_t}$ on every fibre of $\pi'$.

Let
\[\rho':\Theta_{\kb/\Delta} \to \Theta_{\kb}\overset{\rho}{\to} R^1\Phi_*\Theta_{\kx/\kb}\]
be the composition map  which is a map of vector bundles by our assumptions. We claim that $\rho'$ is in fact identically zero:

Let $Z$ be an analytic subspace of $\kb$ of codimension 1 and let
$\ki$ be the corresponding ideal sheaf. Then for any sheaf $\kg$ on $\kx$ there is the base change map
\[ \Phi_\ki:R^1\Phi_*\kg/\ki\cdot R^1\Phi_*\kg \to R^1\Phi_*(\kg\restr{\inverse\Phi Z})\]
which is injective (\cite{Gr-Rem}, Prop. 2, p. 208, p. 209).

For our subspaces $\kb_t$ the naturality of the base change map yields a commutative diagram
\[\xymatrix{  \Theta_{\kb/\Delta}\restr{\kb_t}\ar[r]^-{\rho'}\ar[d]^\isom &\left(R^1\Phi_*\Theta_{\kx/\kb}\right)\restr{\kb_t}\ar@{^(->}[d]\\
 \Theta_{\kb_t}\ar[r]_-{\rho_t} & R^1{\Phi_t}_*\left(\Theta_{\kx_t/\kb_t}\right).}\]

If $t\neq 0$ then we have $\rho_t\equiv 0$ since  $\Phi_t$ is locally trivial. Therefore the map of vector bundles $\rho'\equiv 0$ because it vanishes on the dense open set $\kb\setminus \kb_0$.

Looking again at the diagram for $t=0$ we see that also $\rho_0$ must be zero. Hence $\Phi_0$ is a locally trivial fibration as claimed.\qed

\begin{rem}
By a theorem of Wavrik \cite{wavrik69} condition \refenum{ii} in Proposition \ref{locallytrivial} holds if $h^k(\inverse\Phi(p), \Theta_{\kx_p})$  does not depend on $p\in \kb$ for $k=0$ or $k=1$ and the Kuranishi space of $\inverse\Phi(p)$ is reduced for all $p\in \kb$.  
\end{rem}

\textit{Proof of Theorem \ref{largedeform}.}
Let $\kc$ be a good fibre class and $\pi:\kx\to \Delta$ be a family of compact, complex manifolds. Let $(t_\nu)_{\nu\in \IN}$ be a sequence in $\Delta$ converging to zero such that subspace $V\subsetneq H_1(X,\IQ)$ is $\kc$-fibering on $\kx_{t_\nu}$ for all $\nu$. 

In particular $V$ defines a very good Albanese-Quotient on $\kx_{t}$ for  $t=t_{\nu}$ and hence for all $t\in \Delta$ by Corollary \ref{albaneselimit}. So we may assume that we are in the situation of Lemma \ref{quotientfamily}.

Our assumptions guarantee that we can first apply Proposition \ref{smoothfibration} and then Proposition \ref{locallytrivial}. This concludes the proof of the theorem.\qed

\section{Complex structures on certain Lie-algebras}\label{structureclass}

The aim of this section is to study the possible complex structures on certain types of Lie-algebras. We are particularly interested in the existence of SPTBS in Lie-algebras in view of  Corollary  \ref{stableinvariantdeformation} and Theorem \ref{nillargedeform}.

Using only the dimensions of the subspaces in the descending and ascending central series we will try to give a complete picture for Lie-algebras with commutator subalgebra of dimension at most three. By giving lots of examples we will also show that our classification cannot be improved without considering other properties of the Lie-algebras.
In the case $\dim\kc^1\lg=1$ we can show that there is  a unique complex structure up to isomorphism.

\subsection{Notations and basic results}

In the sequel $\lg$ will denote a nilpotent Lie-algebra and $J$ a complex structure on $\lg$ which is always assumed to be integrable in the sense of Definition \ref{defintegrable}; in particular $\dim \lg$ is always even.

 We continue to use the notation introduced in Section \ref{LICS}, see Section \ref{examsection} for the notation in the examples.

For later reference we collect some basic facts about complex structures and Lie-algebras in the next lemma.
%
%
\begin{lem}\label{basicfacts}
Let $(\lg, J)$ be a real Lie-algebra with a  complex structure.
\begin{enumerate}
\item  Let $V\subset \lg$ be a real subspace. Then the following are equivalent:
\begin{itemize}
\item $V$ is $J$-invariant.
\item $V$ is a complex subspace of $(\lg,J)$ considered as a complex vector space.
\item $ V_\IC=(V_\IC\cap \einsnull\lg)\oplus(V_\IC\cap \nulleins\lg)$.
\end{itemize}
\item If $V\subset \lg$ is $J$-invariant then $\dim_\IR V$ is even. In particular if $V\subset \lg$ is a nontrivial, $J$-invariant subspace contained in a real 2-dimensional subspace $W$ then $V=W$.
\item Let $V$ be a real subspace of $\lg$. Then the largest $J$-invariant subspace of $V$ is
\[\left( (V_\IC\cap \einsnull\lg)\oplus(V_\IC\cap \nulleins\lg)\right)\cap \lg\]
where $\lg$ is identified with $\{x\tensor 1\mid x\in\lg\}\subset \lg_\IC$.
\item For any $x\in \lg$ holds: $ad_x(-)=[x,-]\neq 0 \iff x\notin \kz\lg$.
\item If $\lg$ is 2-step nilpotent then $\kc^1\lg\subset \kz\lg$.
\item Let $V\subset \lg$ be a real subspace with $\dim_\IR V=3$. If $W_1, W_2 \subset V$ are both non-trivial, $J$-invariant subspaces, then $W_1=W_2$.
\item The complex structure is abelian, i.e., $\einsnull \lg$ is an abelian subalgebra of $\lg_\IC$,  if and only if $[x,y]=[Jx,Jy]$ for all $x,y\in \lg$. In this case the ascending central series $(\kz^i\lg)$ is $J$-invariant.
\item If $\kc^1\lg$ contains no $J$-invariant subspace then the complex structure $J$ is abelian.
\end{enumerate}
\end{lem}
\pf We prove only the last three assertions, the rest being clear.
\begin{enumerate}\setcounter{enumi}{5}
\item  Both $W_1$ and $W_2$ have positive, even real dimension, hence 
\[\dim_\IR W_1=\dim_\IR W_2=2\text{ and }W_1\cap W_2\neq 0\] by dimensional reasons. Since the intersection $W_1\cap W_2$ is again $J$-invariant it has also dimension 2 by \refenum{ii} and we have 
\[W_1=W_1\cap W_2=W_2\]
 as claimed.
\item This is a straightforward calculation which can be found for example in \cite{mpps06}.
\item Since we assumed $J$ to be integrable $\einsnull{\lg}$ is a subalgebra of $\lg_\IC$ and hence
\[\kc^1\einsnull \lg\subset \kc^1\lg_\IC\cap \einsnull\lg\subset \kc^1\lg_\IC.\]
If $\kc^1\lg$ contains no $J$-invariant subspace then $\kc^1\lg_\IC\cap \einsnull\lg=0$ by \refenum{i}.
\end{enumerate}
\qed

The main tool in our analysis will be the following technical lemma. 

\begin{lem}\label{keylem}
Let $\lg$ be a nilpotent Lie-algebra with complex structure $J$. Let $\kz_J\lg:= \kz\lg+J\kz\lg$ be the smallest complex subspace of $\lg$ which contains the centre, $\kw_i:= \kz_J\lg \cap \kc^i\lg$ and 
\[\kv^J_{i+1}:= [\kw_i,\lg] = span_\IR\{[w,x]\mid w\in \kw_i, x\in \lg\}.\]
Then:
\begin{enumerate}
\item $\kv^J_{i}$ is a (possibly trivial) $J$-invariant subspace of the $i$-th commutator subalgebra $\kc^i\lg$ and we have $\kv^J_{i+1}\subset \kv^J_{i}$.
\item If there is some $x\in \kw_i\setminus \kz\lg$ then $\kv^j_{i+1}\neq0$.
\item The centre is $J$-invariant if and only if $\kv^J:=\kv^J_1=0$.
\end{enumerate}
\end{lem}
If the complex structure is fixed we will often omit it from the notation.

\pf
We begin with the last assertion: the centre is not fixed by $J$ if and only if for some  $z\in \kz\lg$  the element $Jz$ is not in the centre, which means that there exists some $x\in \lg$ such that $[Jz,x]\neq 0$  which is equivalent to $\kv^J\neq 0$.

For \refenum{i} we only have to show that $x\in \kv^J_i$ implies $Jx\in \kv^J_i$. We will do this on generators of the form  $x=[Jz,y]$ for some $z$ in the centre of $\lg$ and $y\in \lg$. The Nijenhuis tensor then implies
\[ Jx=J[Jz,y]=[Jz,Jy]-[z,y]-J[z,Jy]=[Jz,Jy]\in \kv^J_i.\]

The second assertion follows immediately from the definition of $\kv^J_i$ and \ref{basicfacts} \refenum{iv}.
\qed

In some cases the lemma will enable us to prove that there do not exist complex structures on a certain class of Lie-algebras with the following argument: we assume the existence of a complex structure and then deduce that some odd-dimensional subspace should be invariant under $J$ which is impossible.

We give two applications:

\begin{prop}\label{2nilJnil}
Let $\lg$ be a 2-step nilpotent Lie-algebra. Then every integrable complex structure on $\lg$ is nilpotent. 
\end{prop}
\pf Let $J$ be a complex structure on $\lg$. We will prove the claim by induction on the dimension of $\lg$. We have to show that for some $k$ we have $\kt^k\lg=\lg$ in the minimal torus bundle series or, equivalently, that $\lg/\kt^{k-1}\lg$ is an abelian Lie-algebra.

Since $\lg$ is 2-step nilpotent we have
\[\lg\supset \kz\lg\supset \kc^1\lg \supset0\]
and, by Lemma \ref{keylem}, either $\kv_1$ or $\kz\lg$ is a nontrivial J-invariant subspace of the centre. In particular the centre has dimension at least 2 and if  $\dim \kz\lg=2$  then the centre is $J$-invariant (\ref{basicfacts} \refenum{ii}). This proves the claim if $\lg$ has dimension 4 since the quotient is then necessarily abelian.

In higher dimension we have $\kt^1\lg\neq 0$. If $\lg/\kt^1\lg$ is abelian we are done, else $\lg/\kt^1\lg$ is still 2-step nilpotent and we can use the induction hypothesis.\qed

\begin{prop}\label{Z=C}
Let $\lg$  be a nilpotent Lie-algebra and $k=\dim\kz\lg$. Let $m=k+1$ if $k$ is odd and $m=k+2$ if $k$ is even. Assume that $\kz\lg=\kc^1\lg$ and that one of the following holds:
\begin{enumerate}
\item For every $m$-dimensional subspace $W$, which contains the centre, we have $[W,\lg]=\kc^1\lg$.
\item The map $ad_x:\lg\to \kc^1\lg$ is surjective for all $x\notin \kz\lg$.
\end{enumerate}
If $k=\dim(\kz\lg)$ is not even then there does not exist any integrable complex structure on $\lg$ and if $k=\dim(\kz\lg)$ is even
then any complex structure on such a $\lg$ is nilpotent and 
\[ \lg\supset \kz\lg=\kc^1\lg\supset 0\]
is a SPTBS for $\lg$.
\end{prop}
\pf
Clearly the \refenum{ii} implies the \refenum{i}. Now assume we have a complex structure $J$  on $\lg$ such that the centre is not $J$-invariant. Then $\kz_J\lg:= \kz\lg+J\kz\lg$ is an even dimensional subspace such that $\kz\lg\subsetneq  \kz_J\lg$ and therefore it has dimension at least $m$. The subspace $\kv^J$ is nonempty and by assumption \refenum{i} it is in fact equal to $\kc^1\lg=\kz\lg$. But then  Lemma \ref{keylem} implies that $\kz\lg=\kv^J$ is J-invariant -- a contradiction.

It follows that the centre is $J$-invariant for every complex structure $J$ on $\lg$. Therefore the centre cannot have odd dimension if there exists a complex structure.\qed

In Section \ref{dim3} we will also need the following:
\begin{lem}\label{dim6nn}
 Let $(\lg, J)$ be a nilpotent Lie-algebra with a complex structure, $z\in \kz\lg$ and  $W:=\langle z, Jz\rangle$. Then $ \left(\im ad_{Jz}\right)\cap W=0$.
\end{lem}
\pf Since $\lg$ is nilpotent $Jz\notin \im ad_{Jz}$ and it remains to prove that $z\notin\im ad_{Jz}$. So assume that there is some $x\in \lg$ such that $[Jz,x]=z$. Then with the same calculation as in Lemma \ref{keylem} we have
\[ad_{Jz}(Jx)=[Jz,Jx]=J[Jz,x]=Jz,\]
so $Jz\in \im ad_{Jz}$ which is a contradiction.\qed

\subsection{The case $\dim(\kc^1\lg)=1$}\label{dimcg=1}
Recall that a Heisenberg algebra $H_{2n+1}$ is a nilpotent Lie-algebra which admits a basis $x_1, \dots, x_n, y_1, \dots , y_n, c$ such that $c$ is central and the  structure equations are 
\begin{equation}\label{heis}
\begin{split}
 [x_i,y_i]=-[y_i,x_i]=c, \qquad& i=1,\dots n,\\
[x_i,x_j]=-[y_i,y_j]=0,\qquad &i,j=1,\dots n.
\end{split}
\end{equation}

\begin{prop}\label{propdimcg=1}
Let $\lg$ be a real nilpotent Lie-algebra with $\dim(\kc^1\lg)=1$. Then $\lg$ is 2-step nilpotent and $\lg\isom H_{2n+1}\oplus \IR^m$ is the direct sum  of a Heisenberg algebra with an abelian Lie-algebra with generators $z_1, \dots , z_m$.
Then for any choice of signs
\begin{gather}
Jx_i:=\pm y_i \qquad (i=1,\dots n),\notag\\
 Jz_{2k-1}:=z_{2k}\qquad (k=1,\dots r),\label{heis2}\\
Jc:=z_{2r+1},\qquad J^2=-id_\lg, \notag
\end{gather}
defines a complex structure on $\lg$ and every complex structure on $\lg$ is equivalent to a complex structure of this kind.
Thus, every complex structure is abelian and the ascending central series is a SPTBS for $\lg$. 
\end{prop}
The centre has necessarily even dimension since we assumed $\lg$ to have even dimension.

\pf
Let $c\in \kc^1\lg$ be a generator of the commutator subalgebra. Considering the Lie-bracket as an alternating form on $\lg/\langle c\rangle$   the classification of alternating forms on vector spaces. (\cite{Lang}, XIV, 9) yields the claimed decomposition.

Now assume that $\lg$ admits an integrable complex structure $J$.
The commutator $\kc^1\lg$ cannot contain any complex subspace for dimensional reasons and writing  the Nijenhuis tensor as
\[ [Jx,Jy]-[x,y]=J([Jx,y]+[x,Jy])\]
we see that the left hand side is in $\kc^1\lg$ while the right hand side is not. This yields
\[[Jx,Jy]=[x,y]\]
for all $x,y\in \lg$ and hence any complex structure on $\lg$ is abelian by Lemma \ref{basicfacts} \refenum{vii}.

It is a straightforward calculation to show that $J$ as in \refb{heis2} defines a complex structure and it remains to show that every integrable complex structure can be written in this way with respect to a suitable basis.

Assume that we have $\lg\isom  H_{2n+1}\oplus \IR^{2r+1}$  as above and we are given a complex structure $J$ on $\lg$. The centre $\kz\lg$ is $J$-invariant and we can choose a basis for the centre such that 
\[Jz_{2k-1}=z_{2k}\qquad (k=1,\dots r),\qquad Jc=z_{2r+1} .\]

Let us fix an arbitrary complex subspace $V$ such that $\lg=V\oplus\kz\lg$. The remaining elements of the basis are provided by the following:

\noindent\textbf{Claim:} If $(V,J)$ is a real vector space with a complex structure and
\[[-,-]:V\times V\to \IR\cdot c\] a non-degenerate alternating bilinear form on $V$, such that $[x,y]=[Jx,Jy]$ for all $x,y\in V$ then there exists a basis $x_i, y_i$ of $V$ which satisfies \refb{heis} and \refb{heis2}.

We will prove our claim by induction on the dimension of $V$. First of all we show that we can always find an $x\in V$ such that $[x, Jx] \neq 0$.

We pick a nonzero $a\in V$. If $[a, Ja] \neq0$ then we are done and hence assume $[a, Ja] =0$. There  is some $b\in V$ such that $[a,b]=c$ since the bracket is non-degenerate. If $[b,Jb]\neq 0$ we set $x=b$.

Otherwise  also $[b,Jb]= 0$ and we calculate using $[Ja,Jb]=[a,b]$ 
\begin{align*}
[a+Jb, J(a+Jb)]&=[ a+Jb, Ja-b]\\&=[a,Ja]-[a,b]+[Jb, Ja]-[Jb,b]\\&=-2[a,b]=-2c\neq0, 
\end{align*}
i.e., we can set $x=a+Jb$.

If $[x, Jx] \neq 0$ then $[x, Jx]=\lambda c$ for some $\lambda\in \IR\setminus\{0\}$. Setting $x_1= x/\sqrt{\mid\lambda\mid}$ and  $y_1:=sign(\lambda) Jx_1$ (where $sign(\lambda)$ is the sign of $\lambda$) we have
\[[x_1,y_1]=c \text{ and } Jx_1=sign(\lambda)y_1.\]
If $\dim_\IR V=2$ we are done; otherwise we can apply the induction hypothesis to the subspace orthogonal to $\langle x_1, y_1\rangle$ with respect to the Lie-bracket.
This concludes the proof.\qed

\subsection{The case $\dim(\kc^1\lg)=2$}

We begin with a preliminary lemma.
\begin{lem}\label{CvZ}
Let $\lg $ be a nilpotent Lie-algebra with $\dim(\kc^1\lg)=2$ and $J$ a complex structure on $\lg$. Then at least one of the subspaces $\kz\lg$ and $\kc^1\lg$ is $J$-invariant.
\end{lem}
\pf We have $0\subset \kv\subset \kc^1\lg$ and since $\kv$ is even dimensional and $\dim(\kc^1\lg)=2$, we have either $\kv=0$ or $\kv= \kc^1\lg$. An application of Lemma \ref{keylem} \refenum{iii} concludes the proof. \qed

The classification in this case reads as follows.
\begin{prop}\label{dimcg=2} Let $\lg $ be a nilpotent Lie-algebra with $\dim(\kc^1\lg)=2$. Then every complex structure on $\lg$ is nilpotent and we have the following cases.
\begin{enumerate}
\item If $\lg$ is  3-step nilpotent then the following holds: if one of the subalgebras $\kz^i\lg$ has odd dimension, then $\lg$ does not admit any complex structure. If there exists a complex structure $J$ on $\lg$ then $J$ is abelian and the ascending central series is a SPTBS.
\item If $\lg$ is 2-step nilpotent and $\dim(\kz\lg)$ is odd then 
\[0\subset \kc^1\lg \subset \lg\]
is a SPTBS on $\lg$.
\item If $\lg$ is 2-step nilpotent and $\dim(\kz\lg)$ is even then for every complex structure on $\lg$ either 
\[0\subset \kc^1\lg \subset \lg\quad\text{or}\quad 0\subset \kz\lg \subset \lg\]
is a torus bundle series on $\lg$ but a SPTBS does not exist in general if $\dim(\kz\lg)\geq 4$.
If $\dim(\kz\lg)=2$ the two series coincide and we have a SPTBS.
\end{enumerate}
\end{prop}

\pf
We treat the cases separately:
\begin{enumerate}
\item Assume that we have a complex structure $J$ on $\lg$ and $\lg$ is 3-step nilpotent. It suffices to show that  $\kc^1\lg$ is not a $J$-invariant subspace by Lemma \ref{basicfacts} \refenum{vii}, \refenum{viii} since $\dim_\IR\kc^1\lg=2$. 

Assume the contrary. Writing  $\kc^2\lg=\langle c_2\rangle$ we  then have $\kc^1=\langle c_1:=Jc_2, c_2\rangle$, in particular $\kw_1\neq0$ in the notation of Lemma \ref{keylem}. But $c_2\in \im\left(ad_{c_1}\right)$ and hence, by Lemma \ref{keylem} \refenum{ii}, $\kv_2\neq 0$ which is impossible since $\dim\kc^2\lg=1$.
\item Clearly $\kz\lg$ can never be $J$-invariant for any complex structure on $\lg$ if $\dim \kz\lg$ is odd and hence the assertion follows from Lemma \ref{CvZ}.
\item  We observe that 2-step nilpotency implies $\kc^1\lg\subset\kz\lg$ and hence both subspaces coincide if $\dim\kc^1\lg=\dim\kz\lg=2$. So by Lemma  \ref{CvZ}
\[0\subset \kc^1\lg=\kz\lg \subset \lg\]
is a SPTBS in this case. The remaining assertions follow from Example \ref{ex2ev} given below.
\end{enumerate}
\qed

\begin{exam}\label{ex2ev}
Consider the following 10-dimensional Lie-algebra given by a basis $\lg=\langle e_1, e_2, e_3, f_1, f_2, f_3, z_1, z_2, c_1, c_2\rangle$ and the structure equations 
\[dc^1=e^{12}+f^{12}, \qquad dc^2=e^{13}+f^{13}\]
with respect to the dual basis.
$\lg $ is a  2-step nilpotent Lie-algebra with centre $\kz\lg=\langle  z_1, z_2, c_1, c_2\rangle$ and commutator $\kc^1\lg=\langle   c_1, c_2\rangle$. We give three interable complex structures $J_1$, $J_2$ and $J_3$ on $\lg$ by 
\begin{align*}
J_1e_1&=f_1,&J_2z_1&=e_1,&J_3e_1&=f_1,\\
J_1e_2&=f_2,&J_2z_2&=f_1,&J_3e_2&=f_2,\\
J_1e_3&=f_3,&J_2c_1&=c_2,&J_3e_3&=f_3,\\
J_1z_1&=c_1,&J_2e_2&=e_3,&J_3z_1&=z_2,\\
J_1z_2&=c_2,&J_2f_2&=f_3,&J_3c_1&=c_2.\\
\end{align*}
 The complex structure $J_1$ leaves the commutator invariant but not the centre ($\kv=\kc^1\lg$), $J_2$ the centre but not the commutator while $J_3$ leaves both subspaces invariant. This realises all possible combinations in Proposition \ref{dimcg=2} \refenum{iii}.
\end{exam}

\subsection{The case $\dim(\kc^1\lg)=3$}

\begin{prop}\label{dimcg=3}
Let $\lg$ be a  nilpotent Lie-algebra with $\dim(\kc^1\lg)=3$. Then the following cases can occur:
\begin{enumerate}
\item If $\lg$ is  4-step nilpotent then the following holds: if one of the subalgebras $\kz^i\lg$ has odd dimension, then $\lg$ does not admit any complex structure. If there exists a complex structure $J$ on $\lg$ then $J$ is abelian and the ascending central series is a SPTBS.
\item If $\lg$ is  3-step nilpotent  and $\dim(\kc^2\lg)=\dim(\kc^1\lg\cap \kz\lg)= 1$ then for any complex structure $J$ on $\lg$ we have $\kc^2\lg\cap \kv^J=0$ but $\kv^J\neq 0 $ is possible and there may be nilpotent and non-nilpotent complex structures on the same Lie-algebra.
\item If $\lg$ is 3-step nilpotent such that  $\dim(\kc^1\lg\cap \kz\lg)= 2$ then: 
\begin{enumerate}
\item If $\dim \kz^i\lg$ is even for all $i$ then every complex structure $J$ on $\lg$ is nilpotent and one of the subspaces $\kz\lg$ and $\kc^1\lg\cap \kz\lg$ is $J$-invariant but there is no SPTBS in general.
\item If $\dim \kz^1\lg$ is odd or equal to 2 and $\dim\kz^2\lg$ is even then 
\[ \lg\supset \kz^2\lg\supset\kc^1\lg\cap \kz\lg\supset 0\]
is a SPTBS for $\lg$.
\item If $\dim \kz^2\lg$ is odd there does not exist any complex structure on $\lg$. 
\end{enumerate}
\item 
If $\lg$ is 2-step nilpotent then every complex structure on $\lg$ is nilpotent but in general there does not exist a SPTBS. 
\end{enumerate}
\end{prop}

\pf The last assertion follows immediately from Proposition \ref{2nilJnil} and Example \ref{badex}. The remaining cases \refenum{i}, \refenum{ii} and \refenum{iii} will be treated separately in the following lemmata.\qed

\begin{lem}
Let $\lg$ be a 4-step nilpotent Lie-algebra with $\dim(\kc^1\lg)=3$. If one of the subalgebras $\kz^i\lg$ has odd dimension, then $\lg$ does not admit any complex structure. If there exists a complex structure $J$ on $\lg$ then $J$ is abelian and the ascending central series is a SPTBS.
\end{lem}
\pf Assume that $\lg$ admits a complex structure $J$. 
It suffices to show that  $\kc^1\lg$ does not contain a $J$-invariant subspace by Lemma \ref{basicfacts} \refenum{vii}, \refenum{viii}.

Assume the contrary, i.e., that there is a non-trivial, J-invariant subspace $W\subset \kc^1\lg$ with (necessarily) $\dim_\IR(W)=2$. By Lemma \ref{basicfacts} \refenum{vi} every complex subspace of $\kc^1\lg$ is contained in $W$ and equal to $W$ if it is non-trivial.
Since $\dim \kc^1\lg =3 $ and $\lg$ is 4-step nilpotent we have
\begin{gather*}
\dim \kc^3\lg=1,\qquad \dim \kc^2\lg=2,\\ \kz\lg\cap \kc^1\lg=\kc^3\lg.
\end{gather*}
We will now derive a contradiction.
\begin{description}
\item[Case 1: $W\cap \kc^3\lg=\langle c_3\rangle\neq 0$] Then $Jc_3\in \kc^1\lg\setminus \kz\lg$ and  $\kv_2\neq 0$ by Lemma \ref{keylem} (\textit{ii}). But this in turn implies $\kc^2\lg=\kv_2=W$ for dimensional reasons. Hence $Jc_3\in\kc^2\lg\setminus \kz\lg$ and again by Lemma \ref{keylem} we have $0\neq\kv_3\subset \kc^3\lg$ which is impossible since $\kv_3$ has dimension at least 2.
\item[Case 2: $W\cap \kc^3\lg =0$] 
We construct a basis of $\kc^1\lg $ in the following way: for dimensional reasons $W$ and $\kc^2\lg$ intersect non-trivially and we can choose a nonzero  $c_2\in W\cap \kc^2\lg$. Let $c_1:=Jc_2$ and $c_3\in \kc^3\lg$ be a generator. Then by our assumptions $c_1, c_2, c_3$ is a basis of $\kc^1\lg$. 

There is $f_2\in \lg$ such that $[c_2, f_2]=c_3$. The Nijenhuis tensor together with $Jc_2=-c_1$ yields
\[J(\underbrace{[Jc_2, f_2]+[ c_2, Jf_2]}_{\in \kc^2\lg})= [Jc_2, Jf_2]-[c_2, f_2]\in \kc^2\lg\]
and, since $\kc^2\lg\cap J\kc^2\lg=0$ under our conditions, both sides are equal to zero.
In particular 
\begin{equation}\label{equ1}
[c_1, f_2]=[Jc_2, f_2]=-[ c_2, Jf_2]\in \kc^3\lg
\end{equation}
is central. 
Writing $c_1=[a,b]$ we also have
\begin{equation}\label{equ2}
[c_1, c_2]=[[a,b], c_2]=[ a, [b,c_2]]-[b,[a,c_2]]\in \kc^4\lg=0.
\end{equation}
Since $c_2\in \kc^2\lg\setminus \kz\lg$ we can find $f_1\in \lg$ and $\lambda\in \IR$ such that 
\begin{equation}\label{equ3}
[c_1, f_1]=c_2+\lambda c_3.
\end{equation}
Writing $[f_1,f_2]=\sum_{i=1}^3 \lambda_i c_i$ and using the  Jacobi identity we get
\begin{align*}
 c_3&= [c_2, f_2]= [c_2+\lambda c_3, f_2]\\
&\overset{\text{\refb{equ3}}}{=}[[c_1, f_1], f_2]\\
&= [c_1, [f_1, f_2]] -[f_1, [c_1, f_2]]\\
&\overset{\text{\refb{equ1}}}{=} [c_1, [f_1, f_2]]\\
&=[c_1, \sum_{i=1}^3 \lambda_i c_i]\\
&=\lambda_2 [c_1, c_2]\overset{\text{\refb{equ2}}}{=}0.
\end{align*}
This is a contradiction.
\end{description}

Thus we have shown that $\kc^1\lg$ cannot contain any $J$-invariant subspace which implies in turn that any complex structure will be abelian (Lemma \ref{basicfacts} \refenum{viii}).
\qed
\margincom{Give an example that complex structures can exist?}

\begin{lem}\label{iilem}
Let $\lg$ be a 3-step nilpotent Lie-algebra with $\dim(\kc^1\lg)=3$ and $\dim(\kc^2\lg)=\dim(\kc^1\lg\cap \kz\lg)= 1$. Assume that $J$ is a complex structure on $\lg$. Then $\kc^2\lg\cap \kv^J=0$ but $\kv^J\neq 0 $ is possible and there may be nilpotent and non-nilpotent complex structures on the same Lie-algebra.
\end{lem}

\pf Assume $\kc^2\lg\cap \kv\neq 0$. Then $(J\kz\lg)\cap\kc^1\lg\neq 0$ because $\kv\subset \kc^1\lg$ and  $\kc^2\lg\subset\kz\lg$. Hence $\kv_2\neq 0$ which is impossible since $\dim\kc^2\lg=1$. We prove the second assertion in the following example.\qed

\begin{exam}\label{nilnonnilexam}
Consider the following 10-dimensional Lie-algebra given by a basis $\lg=\langle e_1, e_2, e_3, e_4, e_5, e_6, z, c_1, c_2, c_3\rangle$ and the structure equations 
\begin{gather*}
dc^1=e^{15}+e^{16}+e^{35}+e^{36}\\
dc^2=e^{25}+e^{26}+e^{45}+e^{46}\\
dc^3=e^1\wedge c^1+e^3\wedge c^1+e^2\wedge c^2+e^4\wedge c^2
\end{gather*}
with respect to the dual basis. One can check that $d^2=0$ so the Jacobi-identity holds.
The central filtrations are given by
\begin{gather*}
\lg\supset \kz^2\lg=\langle c_1, c_2, c_3, z \rangle\supset \kz^1\lg=\langle c_3, z\rangle\supset 0\\
\lg\supset \kc^1\lg=\langle c_1, c_2, c_3 \rangle\supset \kc^2\lg=\langle c_3\rangle\supset 0\\
\end{gather*}
and hence we have $\dim \kc^1\lg=3$, $\dim(\kz\lg\cap \kc^1\lg)=\dim \kc^2\lg=1$.
 We give two complex structures $J_1$ and $J_2$  on $\lg$ by 
\begin{align*}
J_1e_1&=e_2,&J_2e_1&=e_2,\\
J_1e_3&=e_4,& J_2e_3&=e_4,\\
J_1e_5&=e_6,&J_2e_5&=z,\\
J_1c_1&=c_2,&J_2c_1&=c_2,\\
J_1z&=c_3,&J_2e_6&=c_3.
\end{align*}
A straightforward calculation shows that both structures are integrable and we have 
\[\kv^{J_1}=0, \qquad \kv^{J_2}=\langle c_1, c_2\rangle\]
 which realises the two possibilities in Lemma \ref{iilem}. Note also that $J_1$ is nilpotent while $J_2$ is not nilpotent.
\end{exam}

\begin{lem}\label{lem32}
Let $\lg$ be a 3-step nilpotent Lie-algebra with $\dim(\kc^1\lg)=3$ and $\dim(\kc^1\lg\cap \kz\lg)= 2$. 
\begin{enumerate}
\item If $\dim \kz^i\lg$ is even for all $i$, then every complex structure $J$ on $\lg$ is nilpotent and one of the subspaces $\kz\lg$ and $\kc^1\lg\cap \kz\lg$ is $J$-invariant  but there is no SPTBS in general.
\item If $\dim \kz^1\lg$ is odd or equal to 2 and $\dim\kz^2\lg$ is even  then 
\[ \lg\supset \kz^2\lg\supset\kc^1\lg\cap \kz\lg \supset 0\]
is a SPTBS for $\lg$.
\item If $\dim \kz^2\lg$ is odd there does not exist any complex structure on $\lg$. 
\end{enumerate}
\end{lem}

\pf
Let $\lg$ be as above and assume that $J$ is a complex structure on $\lg$.
As a first step we show that $\kv\subset\kc^1\lg\cap \kz\lg $. 

Assume the contrary. For dimensional reasons we certainly have 
\[\kv\cap(\kc^1\lg\cap \kz\lg)\neq 0\]
  so that we find some $z\in\kc^1\lg\cap \kz\lg$ such that $Jz\in \kc^1\lg\setminus (\kc^1\lg\cap \kz\lg)$, i.e., $Jz\in \kw_1\setminus \kz\lg$ in the notation of Lemma \ref{keylem} \refenum{ii}. Thus $\kv_2\neq 0$ which is a contradiction if $\dim \kc^2\lg=1$.

Else, if $\dim\kc^2\lg=2$, both $J$-invariant subspaces coincide (Lemma \ref{basicfacts} \refenum{iv}) and 
\[\kv_1= \kv_2=\kc^2\lg= \kc^1\lg\cap \kz\lg\] which contradicts our assumption $\kv\nsubseteq \kc^1\lg\cap \kz\lg$.

Therefore we have $\kv\subset\kc^1\lg\cap \kz\lg $ and at least one of the subspaces $\kz\lg$ and $\kc^1\lg\cap \kz\lg$ will be $J$-invariant. Note that the quotient of $\lg$ by any of these subspaces (which are ideals since they are contained in the centre) is a Lie-algebra with 1-dimensional commutator subalgebra. These have been studied in Section \ref{dimcg=1} and admit only abelian complex structures if the centre is evendimensional and no complex structures if it has odd dimension. This already proves (\textit{iii}) and together with Example \ref{exam322} also (\textit{i}).

If $\dim \kz\lg=2$ the two subspaces coincide and if $\dim \kz\lg$ is odd then $\kz\lg$ can never be invariant whence (\textit{ii}).\qed

In the following two examples we will show that all possibilities in  (\textit{i}) can indeed occur in one Lie-algebra and that we cannot achieve a better result by distinguishing the cases $\dim\kc^2\lg=2$ and $\dim\kc^2\lg=1$

\begin{exam}\label{exam322}
Consider the following 8-dimensional 3-step nilpotent Lie-algebra given by a basis $\lg=\langle e_1, e_2, e_3, c_1, c_2, c_3, z_1, z_2\rangle$ and the structure equations 
\begin{gather*}
dc^1=e^{12}\\
dc^2=e^1\wedge c^1 +e^{23}\\
dc^3=e^2\wedge c^1 -e^{13}
\end{gather*}
with respect to the dual basis. This defines a Lie-algebra since $d^2=0$.
The central filtrations are given by
\begin{gather*}
\lg\supset \kz^2\lg=\langle e_3, c_1, c_2, c_3, z_1, z_2 \rangle\supset \kz^1\lg=\langle c_2, c_3, z_1, z_2\rangle\supset 0\\
\lg\supset \kc^1\lg=\langle c_1, c_2, c_3 \rangle\supset \kc^2\lg=\langle c_2, c_3\rangle\supset 0\\
\end{gather*}
We give three complex structures $J_1$, $J_2$ and $J_3$ on $\lg$ by 
\begin{align*}
J_1e_1&=e_2,&J_2e_1&=e_2,&J_3e_1&=e_2,\\
J_1c_1&=e_3,&J_2c_1&=e_3,&J_3c_2&=c_3,\\
J_1c_3&=c_3,&J_2c_2&=z_1,&J_3z_1&=c_1,\\
J_1z_1&=z_2,&J_2c_3&=z_2,&J_3z_2&=e_3.\\
\end{align*}
It is a straightforward calculation to check that these structures are integrable. The complex structure $J_1$ leaves both the commutator and the centre invariant, $J_2$ the centre but not the commutator,  while $J_3$ leaves the commutator invariant  but not the centre ($\kv_{J_3}=\kc^2\lg$). This realises all possible combinations in Lemma \ref{lem32} (\textit{i}) in a single Lie-algebra.
\end{exam}

\begin{exam}
Consider the following 18-dimensional Lie-algebra given by a basis $\lg=\langle e_1, e_2, e_3, f_1, f_2, f_3, g_1, g_2, g_3, h_1, h_2, h_3, z_1, z_2, c_0, \tilde c_0, c_1, c_2 \rangle$ and the structure equations 
\begin{gather*}
dc^0=g^{12}+g^{23}+h^{12}+h^{23}\\
dc^1=e^{12}+f^{12}\\
\begin{split}
dc^2=&e^{13}+f^{13}\\
&+g^1\wedge c^0+c^0\wedge g^3+h^1\wedge c^0+c^0\wedge h^3\\
&-g^1\wedge \tilde c^0-\tilde c^0\wedge g^3+h^1\wedge \tilde c^0+\tilde c^0\wedge h^3\\
\end{split}
\end{gather*}
with respect to the dual basis.
The central filtrations are given by
\begin{gather*}
\lg\supset \kz^2\lg=\langle e_1, e_2, e_3, f_1, f_2, f_3, c_0, \tilde c_0, c_1, c_2, z_1, z_2 \rangle\supset \kz^1\lg=\langle c_1, c_2, z_1, z_2\rangle\supset 0\\
\lg\supset \kc^1\lg=\langle c_0, c_1, c_2 \rangle\supset \kc^2\lg=\langle c_2\rangle\supset 0
\end{gather*}
and we have $\dim(\kc^1\lg)=3$, $\dim(\kc^1\lg\cap \kz\lg)= 2$ and $\dim(\kc^2\lg)=1$.

We already studied the subalgebra 
\[\gotha=\langle e_1, e_2, e_3, f_1, f_2, f_3, z_1, z_2, c_1, c_2\rangle\]
in Example \ref{ex2ev} and in fact we have an extension
\[ 0\to \gotha\to \lg\to H_7\oplus \IR\to 0\]
where $H_7$ is the 7-dimensional Heisenberg algebra.

On the subspace $V=\langle g_1, g_2, g_3, h_1, h_2, h_3,  c_0, \tilde c_0 \rangle$  complementary to $\gotha$ we give a complex structure $J$ by
\[ Jc_0=\tilde c_0, \qquad Jg_i=h_i, \quad i=1, \dots, 3.\]
Note that for $x,y\in V$ we have 
\[ [x,y]=[Jx,Jy] \text{ and } [x,Jy]=-[Jx,y]\]
and hence the Nijenhuis tensor $[x,y]-[Jx,Jy]+J([Jx,y]+[x,Jy])$ vanishes automatically for these elements in whatever way we choose to extend $J$ to a complex structure on $\lg$. 
Furthermore we have $[V,\gotha]=0$ and therefore any integrable complex structure on $\gotha$ can be combined with $J$ to define an integrable  complex structure on $\lg$.

Combining $J$ with $J_1$,$J_2$ and $J_3$ from Example \ref{ex2ev}  we get three integrable complex structures on $\lg$ which realise the different possibilities in Lemma \ref{lem32} (\textit{i}).
\end{exam}

\section{Applications}\label{applications}

In this section we put together the results obtained so far. In particular we produce a number of classes of nilmanifolds which are closed under deformation in the large and completely describe the situation in real dimension six. 

\subsection{The Main Theorem}
We want to apply   Corollary  \ref{stableinvariantdeformation} and  Theorem \ref{nillargedeform} to our results obtained in Section \ref{structureclass}. The outcome is the following:

\begin{theo}\label{apptheosmall}
Let $M=(\lg, J,\Gamma)$ a nilmanifold with left-invariant complex structure.
Then $M$ and all small  deformations of $M$ are again iterated principal holomorphic torus bundles  coming from a nilmanifold with left-invari\-ant complex structure of type $(\lg,\Gamma)$ if one of the following conditions holds
\begin{itemize}
\item  $\dim \kc^1\lg\leq 2$ .
\item  $\dim \kc^1\lg=3$ and $\lg$ is 4-step nilpotent.
\item $\dim \kc^1\lg=3$, $\lg$ is 3-step nilpotent and $\dim(\kc^1\lg\cap \kz\lg)= 2$.
\end{itemize}
\end{theo}
\pf  We have shown in Proposition \ref{propdimcg=1}, \ref{dimcg=2} and \ref{dimcg=3} respectively that under the above conditions there is a principal torus bundle series for any given complex structure. Combining Theorem \ref{citedolbeault} with Theorem \ref{invariantdeformation} we are done. \qed

In order to proceed to deformations in the large, we need more good fibre classes:
\begin{lem}\label{cg=1goodfibreclass}
If $\lh$ is a nilpotent Lie-algebra with $\dim\kc^1\lh=1$ then  nilmanifolds with left-invariant complex structure of type $(\lh, \Gamma)$ constitute a good fibre class.
\end{lem}
\pf
Recall from Proposition \ref{propdimcg=1} that all complex structures on $\lg$ are equivalent. Hence $h^1(M_J, \Theta_{})=h^{n-1, 1}(M_J)=\dim H^{n-1, 1}(\lg,J)$ is independent of the complex structure $J$. From the same proposition we know that $ 0 \subset \kz\lh\subset\lh$ is a SPTBS, and since complex tori are a good fibre class (see Remark \ref{C-rem}), every deformation in the large is again of the same type by Theorem \ref{nillargedeform}. Thus nilmanifolds with left-invariant complex structure of type $(\lh, \Gamma)$ satisfy both conditions for a good fibre class.\qed

It would be desireable to have more good fibre classes but we can already prove 
\begin{theo}\label{apptheo}
Let $M=(\lg, J,\Gamma)$ a nilmanifold with left-invariant complex structure.
\begin{enumerate}
\item $M$ and all  deformations in the large of $M$ are  iterated principal holomorphic torus bundles coming from a nilmanifold of type $(\lg,\Gamma)$ (with fibre dimensions depending only on $\lg$) if one of the following conditions holds:
\begin{itemize}
\item $\lg$ is abelian.
\item $\dim \kc^1\lg =1$.
\item $\lg$ satisfies the conditions of Proposition \ref{Z=C}.
\item $\dim \kc^1\lg=2$, $\lg$ is 2-step nilpotent and $\dim(\kz\lg)$ is equal to 2 or odd.
\item $\dim \kc^1\lg=2$ and $\lg$ is 3-step nilpotent.
\item $\dim \kc^1\lg=3$, $\lg$ is 3-step nilpotent and furthermore
\begin{itemize}
\item  $\dim(\kc^1\lg\cap \kz\lg)= 2$
\item   $\dim \kz^1\lg$ is odd or equal to 2.
\item   $\dim\kz^2\lg$ is even.
\item  $\dim{\kc^1(\kz^2\lg)}=1$. 
\end{itemize}
The last condition is automatically fulfilled if $\dim \kc^2\lg=1$.
\end{itemize}
\item $M$ and all small  deformations of $M$ are again iterated principal holomorphic torus bundles  coming from a nilmanifold with left-invari\-ant complex structure of type $(\lg,\Gamma)$ if one of the following conditions holds
\begin{itemize}
\item $\lg$  admits a STBS $\lg\supset \ks^1\lg\supset 0$ of length 3.
\item $\lg$ is 3-step nilpotent and admits a STBS of the form $\lg\supset \ks^2\lg\supset\ks^1\lg\supset 0$ such that $\dim(\kc^1(\ks^2\lg))=1$. 
\end{itemize}
\end{enumerate}
\end{theo}
\pf
By Theorem \ref{nillargedeform}, we have to check that in each of the cases there is a SPTBS (resp. STBS) $(\ks^i\lg)_{i=0,\dots, t}$ and that the nilmanifolds of the type $(\ks^{t-1}\lg, \Gamma\cap \exp(\ks^{t-1}\lg))$ constitute a good fibre class. The first assertion either holds by definition or has been proved in Proposition \ref{Z=C}, \ref{propdimcg=1}, \ref{dimcg=2} and \ref{dimcg=3} respectively.

If $t=1$ then $M$ is a torus and the claim follows from \cite{catanese04}. If $t=2$, i.e., $\lg$ admits a SPTBS (or STBS) 
\[\lg\supset \ks^1\lg\supset 0\]
then $M$ is a holomorphic torus bundle over a torus and we are done, since tori constitute a good fibre class (see Remark \ref{C-rem}). This is always the case if $\dim \kc^1\lg\leq 2$ and $\lg $ is not 3-step nilpotent and if $\lg$ satisfies the conditions of Proposition \ref{Z=C}.

If $\lg$ is 3-step nilpotent and has a SPTBS $(\ks^i\lg)_{i=0,\dots, 3}$, then we have to show that the nilmanifolds with Lie-algebra $\ks^2\lg$ constitute a good  fibre class. 

By Lemma \ref{cg=1goodfibreclass} it suffices to check that $\dim{\kc^1(\ks^2\lg)}=1$ in the remaining cases. For a 3-step nilpotent Lie-algebra we have
\[\kc^1(\ks^2\lg)\subset\kc^1(\kz^2\lg)\subset\kc^2\lg\]
and are done if $\dim \kc^2\lg=1$. In the remaining cases the property holds by assumption. This concludes the proof.\qed
\margincom{The theorem should hold also in the  remaining case "$\dim \kc^1\lg=3$ and $\lg$ is 4-step nilpotent" resp. for all $\lg$ that admit only abelian cx structures!}

\subsection{Deformations and geometric structure in dimension 3}\label{dim3}

In this section we give a fairly complete classification of the geometric types of nilmanifolds with left-invariant complex structure in complex dimension 3 and determine their deformations. 

In \cite{magnin86} Magnin gave a classification of real nilpotent Lie-algebra in dimension at most 7 and in particular showed that in real  dimension 6 there exist only 34 different isomorphism types while in higher dimension there are always continuous families.

Salamon showed in \cite{salamon01} that only 18 of these  6-dimensional real nilpotent Lie-algebras admit a complex structure and Ugarte studied in detail the possible nilpotent and abelian structures (\cite{math.DG/0411254}, Theorem. 2.9); a part of these results has been reproved in Section \ref{structureclass}. Following Ugarte's notation we give the list of  6-dimensional real Lie-algebras admitting complex structures:
\begin{align*}
\lh_1 &= (0,0,0,0,0,0), & \lh_{10} &= (0,0,0,12,13,14),\\
\lh_2 &= (0,0,0,0,12,34), & \lh_{11} &= (0,0,0,12,13,14+23),\\
\lh_3 &= (0,0,0,0,0,12+34), & \lh_{12} &= (0,0,0,12,13,24),\\
\lh_4 &= (0,0,0,0,12,14+23), & \lh_{13} &= (0,0,0,12,13+14,24),\\
\lh_5 &= (0,0,0,0,13+42,14+23), & \lh_{14} &= (0,0,0,12,14,13+42),\\
\lh_6 &= (0,0,0,0,12,13), & \lh_{15} &= (0,0,0,12,13+42,14+23),\\
\lh_7 &= (0,0,0,12,13,23), & \lh_{16} &= (0,0,0,12,14,24),\\
\lh_8 &= (0,0,0,0,0,12), & \lh_{19}^- &= (0,0,0,12,23,14-35),\\
\lh_9 &= (0,0,0,0,12,14+25), & \lh_{26}^+ &= (0,0,12,13,23, 14+25).
\end{align*}

Our classification then reads as follows:
\begin{theo}\label{3dim-class}
Let $M_J=(\lg , J, \Gamma)$ be a complex 3-dimensional nilmanifold with left-invariant complex structure. 
\begin{enumerate}
 \item If $\lg$ is not in $\{\lh_7, \lh_{19}^-, \lh_{26}^+\}$, then $\lg$ admits a SPTBS and hence $M_J$ has the structure of an iterated principal holomorphic torus bundle. We list the possibilities in the following table:
\begin{center}
\begin{tabular}{c|c|c}
 base& fibre & corresponding Lie-algebras\\
\hline
3-torus & - & $\lh_1$\\
2-torus & elliptic curve & $\lh_2, \lh_3  , \lh_4, \lh_5, \lh_6$\\
elliptic curve & 2-torus & $\lh_8$\\
Kodaira surface & elliptic curve & $\lh_9, \lh_{10}, \lh_{11}, \lh_{12}, \lh_{13}, \lh_{14}, \lh_{15}, \lh_{16}$\\
\end{tabular}
\end{center}
Every deformation in the large is of the same type.
\item 
If   $\lg=\lh_{19}^-$ or $\lg=\lh_{26}^+$  then $\lg$ admits a STBS and $M_J$ can be described as 2-torus bundle over an elliptic curve but there is no principal torus bundle structure. Every deformation in the large is of the same type.
\item If $\lg=\lh_7$ then there is a dense subset of the space of all left-invariant complex structures for which $M$ admits the structure of principal holomorphic bundle of elliptic curves over a Kodaira surface but this is not true for all complex structures.
\end{enumerate}
\end{theo}

\pf First of all note that every nilpotent Lie-algebra of real dimension at most 4 which admits complex structures gives rise to a good fibre class, since the only possibilities are elliptic curves, 2-dimensional tori and Kodaira surfaces. So for \refenum{i} and \refenum{ii} it suffices to exhibit SPTBS resp. STBS with the appropriate dimensions. 

We start with the cases in \refenum{i}.  The calculation of the dimensions of the descending and ascending series is straightforward and we find that all cases are covered by the results in Section \ref{structureclass}:
\begin{center}
\begin{tabular}{l|c|c|c|c|c}
 Proposition & \ref{dimcg=1} &    \ref{dimcg=2} (\textit{i}) & \ref{dimcg=2} (\textit{ii})&  \ref{dimcg=2} (\textit{iii})&\ref{dimcg=3} (\textit{iii}) (b) \\
\hline
Lie-algebra &
$\lh_3$, $\lh_8$ & $\lh_9$ &  $\lh_6$&  $\lh_2$, $\lh_4$, $\lh_5$ &  \begin{minipage}[t]{3cm}
$\lh_{10}$, $\lh_{11}$, $\lh_{12}$, $\lh_{13}$,\\
 $\lh_{14}$,  $\lh_{15}$, $\lh_{16}$ 
\end{minipage}
\end{tabular}
\end{center}
The Lie-algebra $\lh_1$ corresponds to the complex 3-dimensional torus.

The Lie-algebras $\lh_{19}^-$ and $\lh_{26}^+$ do not admit any nilpotent complex structure (see \cite{salamon01} or note that the centre has real dimension 1) and therefore a corresponding nilmanifold can never admit an iterated principal holomorphic torus bundle structure. Thus it remains to exhibit a STBS in both cases.
\begin{description}
\item[The case  $\lh_{26}^+$] We claim that $\lg\supset \kc^1\lg:=\langle e_3, \dots, e_6\rangle\supset 0$ is a stable torus bundle series. Since the commutator is always a rational subspace with respect to every rational structure, we only need prove that for every complex structure $J$ on $\lh_{26}^+$ the subspace $\kc^1\lg$ is $J$-invariant. 

Let $Je_6=\sum{i=1}^6\lambda_i e_i$, $W:=\kz_J\lh_{26}^+=\langle e_6, Je_6\rangle$ and, as in Lemma \ref{keylem}, $\kv:=\im ad_{Je_6}$. By Lemma \ref{dim6nn} we know that $\kv\cap W=0$ and, applying $ad_{Je_6}$ to $e_4$ and $e_5$, we see that $\lambda_1=\lambda_2=0$. As $W\subset \kc^1\lg$, $\kv \subset \kc^1\lg$ and their intersection is $0$ so $\kc^1\lg=\kv\oplus W$ is $J$-invariant.
\item[The case  $\lh_{19}^-$] We claim that $\lg\supset V:=\langle e_2, e_4, e_5, e_6\rangle\supset 0$ is a stable torus bundle series. With the same proof as in the first case (replacing $e_2$ by $e_3$) we see that $V$ is $J$-invariant for every complex structure in $\lh_{19}^-$.

It is straightforward to check that
\[V=\{x\in \lh_{19}^-\mid (ad_x)^2=0\}\]
which implies that $V$ is rational with respect to every rational structure.\margincom{elaborate this!}
\end{description}

It remains to prove \refenum{iii}.
A nilmanifold of type $M_J=(\lh_7,J, \Gamma)$ admits a structure of iterated principal holomorphic torus bundle if and only if the subspaces $\kt^1\lh_7$ and $\kt^2\lh_7$ are $\Gamma$-rational. This is clearly the case if the complex structure is rational but not always as seen in Example \ref{badex}. A quotient of $\lh_7$ by a 2-dimensional subspace of the centre can not be abelian, so $\kt^2\lh_7\neq\lh_7$ and, if there is any iterated bundle structure, we have a bundle of elliptic curves over  a Kodaira surface.\qed

\begin{rem}
The real Lie-algebra underlying the Iwasawa manifold is isomorphic to $\lh_5$ and hence we have in particular proved that every deformation in the large of the Iwasawa manifold is a nilmanifold with left-invariant complex structure. It would be interesting to know what are the deformations in the large of nilmanifolds with Lie-algebra $\lh_7$.
\end{rem}

%

\end{document}